\documentclass[preprint,12pt]{elsarticle}

\usepackage{amssymb}
\usepackage{amsthm}
\usepackage{amsmath}

\usepackage{comment}

\usepackage{xcolor}

\def\FigureScale{0.85}

\usepackage[section]{placeins}

\usepackage{tabularx}

\usepackage{eurosym}

\journal{International Journal of Electrical Power \& Energy Systems}

\begin{document}

\begin{frontmatter}

%% Title, authors and addresses

%% use the tnoteref command within \title for footnotes;
%% use the tnotetext command for theassociated footnote;
%% use the fnref command within \author or \address for footnotes;
%% use the fntext command for theassociated footnote;
%% use the corref command within \author for corresponding author footnotes;
%% use the cortext command for theassociated footnote;
%% use the ead command for the email address,
%% and the form \ead[url] for the home page:
%% \title{Title\tnoteref{label1}}
%% \tnotetext[label1]{}
%% \author{Name\corref{cor1}\fnref{label2}}
%% \ead{email address}
%% \ead[url]{home page}
%% \fntext[label2]{}
%% \cortext[cor1]{}
%% \affiliation{organization={},
%%             addressline={},
%%             city={},
%%             postcode={},
%%             state={},
%%             country={}}
%% \fntext[label3]{}

\title{Impact of large-scale hydrogen electrification and retrofitting of natural gas infrastructure on the European power system}

%% use optional labels to link authors explicitly to addresses:
%% \author[label1,label2]{}
%% \affiliation[label1]{organization={},
%%             addressline={},
%%             city={},
%%             postcode={},
%%             state={},
%%             country={}}
%%
%% \affiliation[label2]{organization={},
%%             addressline={},
%%             city={},
%%             postcode={},
%%             state={},
%%             country={}}

\author[inst1]{Germán Morales-España}
\author[inst3]{Ricardo Hernández-Serna}
\author[inst1,inst2]{Diego A. Tejada-Arango}
\author[inst1]{Marcel Weeda}

\affiliation[inst1]{organization={TNO Energy \& Materials Transition},%Department and Organization
            addressline={Radarweg 60}, 
            city={Amsterdam},
            postcode={1043 NT}, 
            country={The Netherlands}}

\affiliation[inst2]{organization={Universidad Pontificia Comillas},%Department and Organization
            addressline={C. Alberto Aguilera, 23}, 
            city={Madrid},
            postcode={28015}, 
            country={Spain}}

\affiliation[inst3]{organization={Eurus Energy Europe BV},%Department and Organization
            addressline={Parnassusweg 821b}, 
            city={Amsterdam},
            postcode={1082 LZ}, 
            country={The Netherlands}}

\begin{abstract}
%% Text of abstract
In this paper, we aim to analyse the impact of hydrogen production decarbonisation and electrification scenarios on the infrastructure development, generation mix, $CO_{2}$ emissions, and system costs of the European power system, considering the retrofit of the natural gas infrastructure. We define a reference scenario for the European power system in 2050 and use scenario variants to obtain additional insights by breaking down the effects of different assumptions. The scenarios were analysed using the European electricity market model COMPETES, including a proposed formulation to consider retrofitting existing natural gas networks to transport hydrogen instead of methane. According to the results, 60\% of the EU's hydrogen demand is electrified, and approximately 30\% of the total electricity demand will be to cover that hydrogen demand. The primary source of this electricity would be non-polluting technologies. Moreover, hydrogen flexibility significantly increases variable renewable energy investment and production, and reduces $CO_{2}$ emissions. In contrast, relying on only electricity transmission increases costs and $CO_{2}$ emissions, emphasising the importance of investing in an $H_{2}$ network through retrofitting or new pipelines. In conclusion, this paper shows that electrifying hydrogen is necessary and cost-effective to achieve the EU's objective of reducing long-term emissions.

\end{abstract}

%%Graphical abstract
\begin{graphicalabstract}
\includegraphics{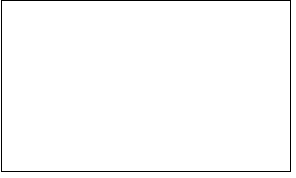}
\end{graphicalabstract}

%%Research highlights
\begin{highlights}
\item A new formulation is proposed for natural gas network retrofitting to hydrogen ($H_2$) 
\item 30\% of EU electricity in 2050 will supply 60\% of the $H_2$ needs
\item Flexible electricity-based $H_2$ production boosts renewables by 48\%
\item Using electricity for $H_2$ production cuts $CO_2$ emissions by 35\%
\item 11\% gas network retrofit covers $H_2$ demand with no new pipelines required
\end{highlights}

\begin{keyword}
%% keywords here, in the form: keyword \sep keyword
Hydrogen \sep Electrification \sep Retrofit \sep Steam Methane Reforming 
%% PACS codes here, in the form: \PACS code \sep code
%%\PACS 0000 \sep 1111
%% MSC codes here, in the form: \MSC code \sep code
%% or \MSC[2008] code \sep code (2000 is the default)
\MSC 90C05 %%\sep 1111
\end{keyword}

\end{frontmatter}

%% \linenumbers

%% main text
%% For citations use: 
%%       \citet{<label>} ==> Jones et al. [21]
%%       \citep{<label>} ==> [21]
%%

%---------------------------------------------------------------------------
%                Introduction
%---------------------------------------------------------------------------

\section{Introduction}
\label{sec:intro} 
\subsection{Background}
The European Union (EU) aims to become carbon-neutral by 2050. This goal is at the heart of the European Green Deal and aligns with the EU's commitment to increase global climate action according to the Paris Agreement. The electrification of end-use services in the transport, residential, and industrial sectors--coupled with the decarbonisation of electricity generation--is one of the essential options for achieving $CO_{2}$ emission reduction targets and climate change mitigation \citep{Williams2012Technology}. The transport and residential sectors can be directly coupled to the power system by adopting electric end-use technologies, such as heat pumps in the residential sector and electric vehicles in the transport sector. Nevertheless, various energy vectors will likely play a role in decarbonising different sectors in the net-zero future. One such vector is green hydrogen, produced using renewable energy \citep{Parra2019review}. For instance, low-carbon hydrogen has been identified as a valuable energy vector for end uses where it is one of the most efficient solutions in decarbonisation, or there is no option for direct electrification, i.e., hydrogen-intensive industry, high-temperature processes, and long-distance heavy transport \citep{Napp2014review}. Furthermore, hydrogen ($H_{2}$) can be an essential long-term energy storage option in 100\% renewable power systems \citep{Le2023}.

\subsection{Electrification and hydrogen}
This emerging electrification trend across industrial processes, electric vehicles, heat pumps, and green $H_2$ production will place additional demands on the power system, requiring significant changes in its planning and operation. Previous research has demonstrated the impact on power systems of electrifying new sectors. For instance, electrification presents opportunities for flexibility, such as the bi-directional charging of electric vehicles, which can help facilitate the integration of renewable energy sources \citep{ramirez2016}. In addition, analyses by \citet{Taljegard2019Impacts} for the Scandinavian countries and Germany, and \citet{Loschan2023EVFlexibility} for Austria have found that electrification of the transportation sector would lead to increased electricity demand, which would be met mainly by wind and thermal plants. Moreover, \citet{Gryparis2020Electricity} found that electrification of vehicles will support efforts to reduce carbon emissions--but not as fast as expected, as a significant percentage of electricity generation in the EU is still based on fossil fuels. This shows the generation mix is an essential consideration to guarantee electrification leads to an actual reduction in $CO_2$ emissions. Although electrifying on a large scale can aid in reducing the carbon footprint of the power system up to only a certain point, it remains a crucial component of the strategy, especially with $H_2$ incorporation. For example, \citet{Sasanpour2021Strategic} found that compared to European power systems that do not use $H_2$, EU Countries that incorporate $H_2$ can obtain a 14-16\% reduction in their total system costs. Moreover, \citet{Pietzcker2021Tightening} found that electrification would help the EU power system meet its decarbonization goals, while  \citet{Lux2020supply} and also \citet{Moser2020sensitivity} found that electrolysis provides flexibility to the power system. These studies suggest that electrification is required for the EU power system to meet its decarbonization goals. However, previous analyses on EU electrification did not fully consider the combined impact of demand response, heating/transport electrification, $H_2$ decarbonisation, and retrofitting gas networks for exclusively $H_2$ transport. This paper considers all of these elements, providing a combined picture with detailed insights into the EU electricity sector of 2050.

\subsection{Modelling the retrofitting of natural gas networks}
In light of the growing need for efficient and sustainable $H_2$ transportation discussed in the previous section, retrofitting natural gas networks for exclusively $H_2$ transport presents a viable solution \citep{Sandana2022}. Therefore, this option should be considered in the optimisation model when analysing the large-scale electrification of $H_2$ production. Nevertheless, some of the state-of-the-art energy planning models in the literature, such as SpineOpt \citep{Ihlemann2022}, TIMES \citep{JRC117820}, and COMPETES \citep{sijm2017a} do not consider this option in their investment decisions. PyPSA-eu \citep{HORSCH2018207} recently added the retrofit as an option in their model \citep{NEUMANN2023}; however, it does not consider the different levels of retrofit that can be developed in a pipeline. Retrofitting the natural gas networks mainly leads to different levels of retrofit, from compressor upgrades to pipeline reinforcements \citep{CRISTELLO2023}. Considering different levels of retrofitting costs in the optimisation model is not straightforward. One option is to formulate a set of constraints considering binary variables; however, this option leads to a Mixed-integer programming (MIP) problem, which is harder to solve in large-scale optimisations such as the European power system. Hence, there is a need for a Linear Programming (LP) formulation that includes different levels of retrofitting as an investment option. To this end, this paper proposes an LP mathematical formulation to incorporate this option, which can be adapted to any energy planning model.

\subsection{Contribution}
\label{subsec:contribution} 
This paper addresses two gaps in the existing literature: a lack of comprehensive analysis of electrification scenarios with $H_2$ decarbonisation and no linear mathematical formulation for different levels of retrofitting natural gas networks for $H_2$ transport. Our contribution includes a detailed power system analysis for the EU and selected countries, measuring the impact of electrification on infrastructure development, generation mix, $CO_{2}$ emissions, and power system costs. We developed a reference scenario for 2050 and analysed it with an optimisation model called COMPETES, including a proposed formulation for retrofitting natural gas infrastructure for $H_2$ transport. Therefore, this paper aims to answer the following research questions:

\begin{itemize}
    \item Is it possible to effectively include the retrofit of natural gas networks as an investment option in energy planning optimisation models?
    \item What are the impacts on the total costs and $CO_2$ emissions of retrofitting the existing gas infrastructure for $H_2$ transport in the EU by 2050?
    \item How do investments in new electrolysers, hydrogen transmission, and storage infrastructure impact total $CO_2$ emissions compared to scenarios where these investments are not made?
\end{itemize}

These questions have been comprehensively analysed and answered throughout this paper, and the results have yielded valuable insights discussed in the conclusion section \ref{sec:conclusions}.
\FloatBarrier
%---------------------------------------------------------------------------
%                Method
%---------------------------------------------------------------------------
\section{Method}
\label{sec:method} 

\begin{figure}[h]
    \centering
    \includegraphics[width=\FigureScale\textwidth]{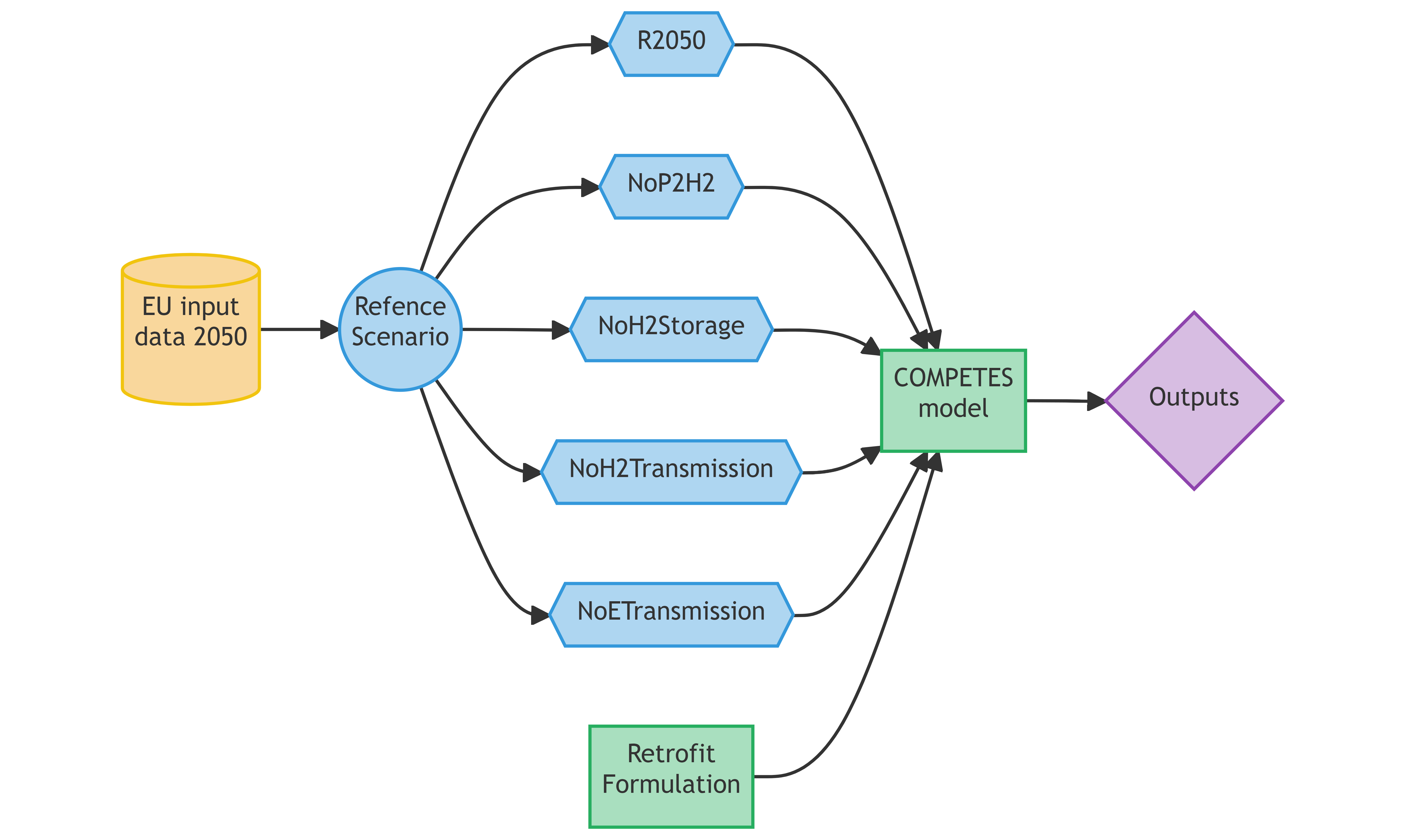}
    \caption{Methodology overview}
    \label{fig:flow_chart}
\end{figure}

This paper studies the impact of large-scale electrification of $H_2$ production in the European electricity sector using a model-based analysis designed to quantify the effect on infrastructure development, generation mix, $CO_{2}$ emissions, and power system costs. Figure \ref{fig:flow_chart} shows an overview of the methodology used to achieve this purpose. First, the input data at the European level helped to define a reference scenario for 2050 (R2050); see Section \ref{sec:refer_sc}. In addition, scenario variants (NoP2H2, NoH2Storage, NoH2Transmission, and NoETransmission) were analysed to obtain additional insights by breaking down the effect of different assumptions in the reference scenario. Section \ref{subsec:H2electr_ScnVar} describes the scenario variants in more detail. Then, all the scenarios were run in the European electricity market model (COMPETES), enhanced with a novel formulation considering the retrofitting of existing natural gas networks to carry $H_{2}$ instead of Methane. The proposed formulation for the retrofit modelling is shown in Section \ref{sec:retrofit_formulation}, while the main COMPETES features are described in Section \ref{sec:opt_model}. Finally, the most relevant outputs of the model are shown in Section \ref{sec:H2electr}, which also discusses the impact of $H_2$ electrification in COMPETES for the reference scenario and its variants.

\FloatBarrier
\subsection{Formulation to retrofit the natural gas network to transport only $H_{2}$}
\label{sec:retrofit_formulation}

One viable techno-economical option to enable the future transport of $H_2$ is retrofitting the existing natural gas networks to transport only $H_{2}$. Therefore, incorporating this option in the optimisation models is vital in analysing the large-scale electrification of $H_2$ production. This section shows the mathematical formulations that can be included in this option in electricity market models such as COMPETES. It is worth noting that the proposed formulation is model agnostic and can be adapted to any other energy planning model.

The retrofit considers a piecewise linear cost curve to account for the investments when increasing the $H_2$ transport capacity in the existing network in the optimisation model. In addition, the retrofit formulation allows accounting for different levels of $H_2$ compression, which increase the $H_2$ energy transport capacity of the network. The mathematical formulation of the retrofitting modelling is as follows:

\begin{align}
\overline{p}_{l}^{1\mathrm{H2}} & \leq\eta^{1}\overline{P}_{l}^{\mathrm{CH4}}\qquad\forall l\label{eq:1st-retrofit}\\
\overline{p}_{l}^{2\mathrm{H2}} & \leq\left(\eta^{2}-\eta^{1}\right)\frac{\overline{p}_{l}^{1\mathrm{H2}}}{\eta^{1}}\qquad\forall l\label{eq:2nd-retrofit}\\
p_{lt}^{\mathrm{H2}} & \leq\left(\overline{P}_{l}^{\mathrm{H2}}+\overline{p}_{l}^{1\mathrm{H2}}+\overline{p}_{l}^{2\mathrm{H2}}+\overline{p}_{l}^{\mathrm{H2}}\right)\Delta t\qquad\forall l,t\label{eq:H2-capacity}
\end{align}

Where $t$ and $l$ are indices for time periods and pipelines, respectively. Parameters $\overline{P}_{l}^{\mathrm{H2}}$ and $\overline{P}_{l}^{\mathrm{CH4}}$ are the initial installed capacities for hydrogen ($H_{2}$) and natural gas (methane, CH4), respectively; $\Delta t$ is the time duration; and $\eta^{1}$ and $\eta^{2}$ are the efficiencies for the first and second $H_{2}$ retrofit, respectively. By definition, $\eta^{2}>\eta^{1}$ since the second retrofit includes higher compression, resulting in higher energy content.

The variable $\overline{p}_{l}^{1\mathrm{H2}}$ is the capacity of the first retrofit which is limited by the initial installed capacity of the natural gas pipeline $\overline{P}_{l}^{\mathrm{CH4}}$ (\ref{eq:1st-retrofit}). The variable $\overline{p}_{l}^{2\mathrm{H2}}$ is the extra capacity added to $\overline{p}_{l}^{1\mathrm{H2}}$ due to higher compression and is limited by (\ref{eq:2nd-retrofit}). Notice that if the natural gas pipeline is completely re-purposed with the second retrofit (higher
compression), $\overline{p}_{l}^{1\mathrm{H2}}+\overline{p}_{l}^{2\mathrm{H2}}=\eta^{2}\overline{P}_{l}^{\mathrm{CH4}}$.

Constraint (\ref{eq:H2-capacity}) limits the $H_{2}$ flow variable $p_{lt}^{\mathrm{H2}}$ to the initial $H_{2}$ capacity $\overline{P}_{l}^{\mathrm{H2}}$, plus the first and second retrofits $\overline{p}_{l}^{1\mathrm{H2}}+\overline{p}_{l}^{2\mathrm{H2}}$, plus the new $H_{2}$ capacity investment variable $\overline{p}_{l}^{\mathrm{H2}}$.

The retrofitted $H_2$ transport capacity decreases the natural gas transport capacity:

\begin{align}
p_{lt}^{\mathrm{CH4}} & \leq\left(\overline{P}_{l}^{\mathrm{CH4}}-\frac{\overline{p}_{l}^{1\mathrm{H2}}}{\eta^{1}}\right)\Delta t\qquad\forall l,t\label{eq:CH4-capacity}
\end{align}

Where $p_{lt}^{CH4}$ is the variable for the natural gas flow. If the first retrofit takes place at its maximum potential, $\overline{p}_{l}^{1\mathrm{H2}}=\eta^{1}\overline{P}_{l}^{\mathrm{CH4}}$ from (\ref{eq:1st-retrofit}), then (\ref{eq:CH4-capacity}) enforces that the maximum available flow capacity for natural gas is zero. Notice that (\ref{eq:1st-retrofit}) is not affected by the second retrofit since the first retrofit uses the initial natural gas pipeline capacity, and the second retrofit uses the same capacity but with higher energy content due to the higher $H_{2}$ compression. The non-negative constraints for all variables are also included.

\begin{align}
\overline{p}_{l}^{1\mathrm{H2}},\overline{p}_{l}^{2\mathrm{H2}},\overline{p}_{l}^{\mathrm{H2}} & \geq0\label{eq:nonnegative-invstment}\\
p_{lt}^{\mathrm{H2}},p_{lt}^{CH4} & \geq0\label{eq:nonnegative-H2flow}
\end{align}

The hydrogen $p_{lt}^{\mathrm{H2}}$ and methane $p_{lt}^{\mathrm{CH4}}$ flows appear in their nodal energy balances \citep{Koirala2021}. The retrofit  $\overline{p}_{l}^{1\mathrm{H2}}+\overline{p}_{l}^{2\mathrm{H2}}$
and investment $\overline{p}_{l}^{\mathrm{H2}}$ variables appear in the objective function with their respective annualized capital expenditure (CAPEX), where the total retrofit and investment cost $c^{\mathrm{TotH2}}$ is given by

\begin{align}
c^{\mathrm{TotH2}} & =\sum_{l}\left(\overline{p}_{l}^{1\mathrm{H2}}C_{l}^{1\mathrm{H2}}+\overline{p}_{l}^{2\mathrm{H2}}\Delta C_{l}^{2\mathrm{H2}}+\overline{p}_{l}^{\mathrm{H2}}C_{l}^{\mathrm{H2}}\right)\\
\Delta C_{l}^{2\mathrm{H2}} & =\frac{\eta^{2}C_{l}^{2\mathrm{H2}}-\eta^{1}C_{l}^{1\mathrm{H2}}}{\eta^{2}-\eta^{1}}\qquad\forall l\label{eq:Capex-Delta2retrofit}
\end{align}

\begin{figure}[h]
    \centering
    \includegraphics[width=\FigureScale\textwidth]{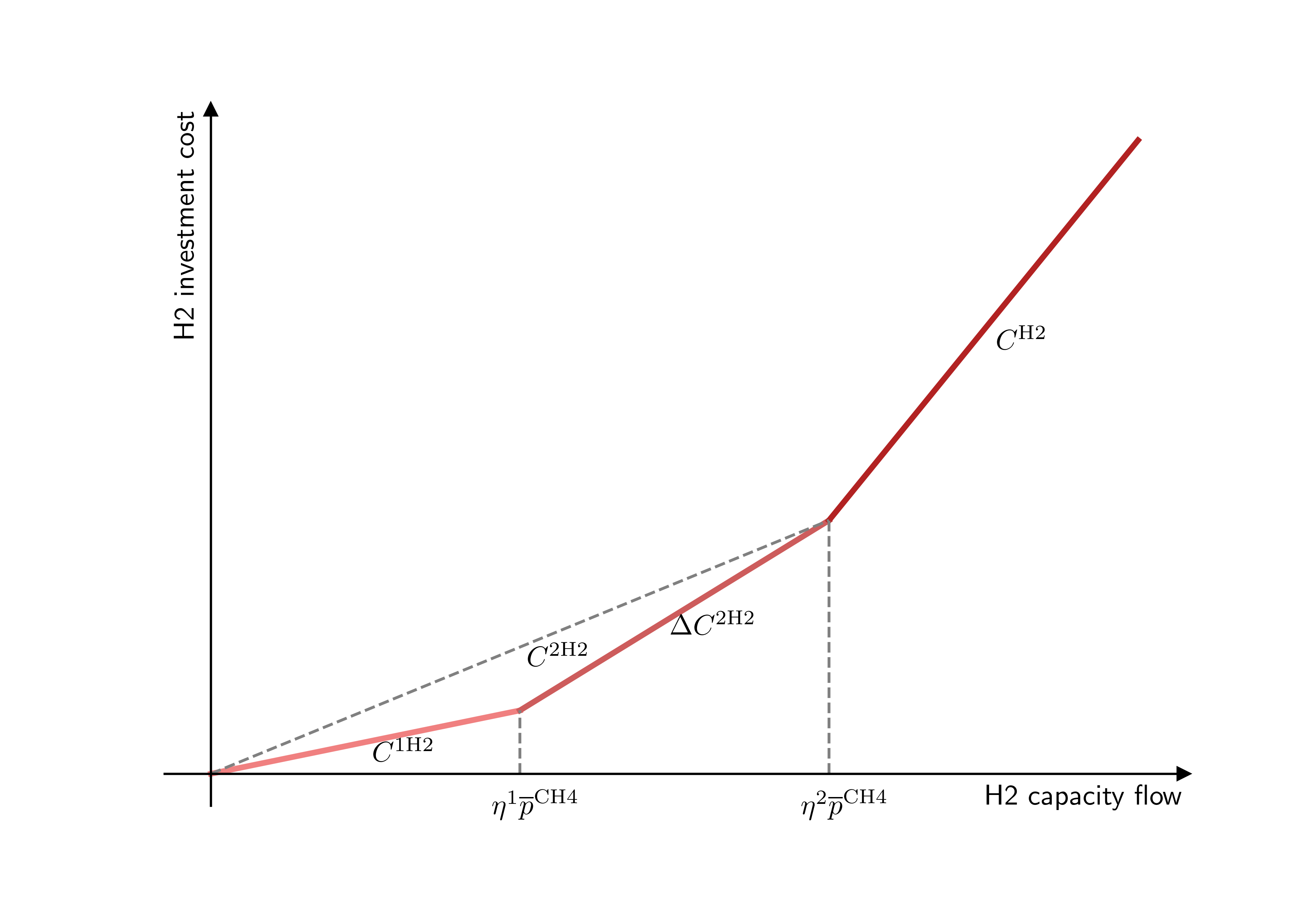}
    \caption{$H_{2}$ retrofit and investment cost}
    \label{fig:H2_retrofit_investment_cost}
\end{figure}

Where $C_{l}^{1\mathrm{H2}}$, $C_{l}^{2\mathrm{H2}}$, and $C_{l}^{\mathrm{H2}}$
are the annualized Capex for the first-retrofit, second-retrofit, and new-pipeline $H_{2}$ investments, respectively. The parameter $\Delta C_{l}^{2\mathrm{H2}}$ is the extra annualized investment cost to reach the second retrofit and is defined in (\ref{eq:Capex-Delta2retrofit}). Figure~\ref{fig:H2_retrofit_investment_cost} shows the relationship between all these annualized investment costs.

The main advantage of this modelling proposal is that it effectively determines the optimal retrofitting decision while keeping the formulation linear. Therefore, the computational burden for a large-scale optimisation model is lower than for formulations using MIP approaches. 

Finally, it is essential to highlight that the natural gas network usually includes multiple parallel pipelines connecting countries. Therefore, the results of this modelling proposal can be interpreted as the number of pipelines that require retrofitting for its implementation.
\FloatBarrier

%---------------------------------------------------------------------------
%                Optimisation Model
%---------------------------------------------------------------------------
\subsection{Optimisation model description}
\label{sec:opt_model} 
In order to address the research questions in Section~\ref{sec:intro}, this paper uses the optimisation model COMPETES, which is a power system optimisation and economic dispatch model that seeks to meet European electricity demand at minimum social costs (i.e., maximising social welfare) within a set of techno-economic constraints – including policy targets/restrictions – of generation units and transmission interconnections across European countries and regions.

COMPETES solves a transmission and generation capacity expansion problem to determine and analyse least-cost capacity investments with perfect competition formulated as a linear programme to optimise the system's generation capacity additions and economic dispatch.

\begin{figure}[h]
    \centering
    \includegraphics[width=\FigureScale\textwidth]{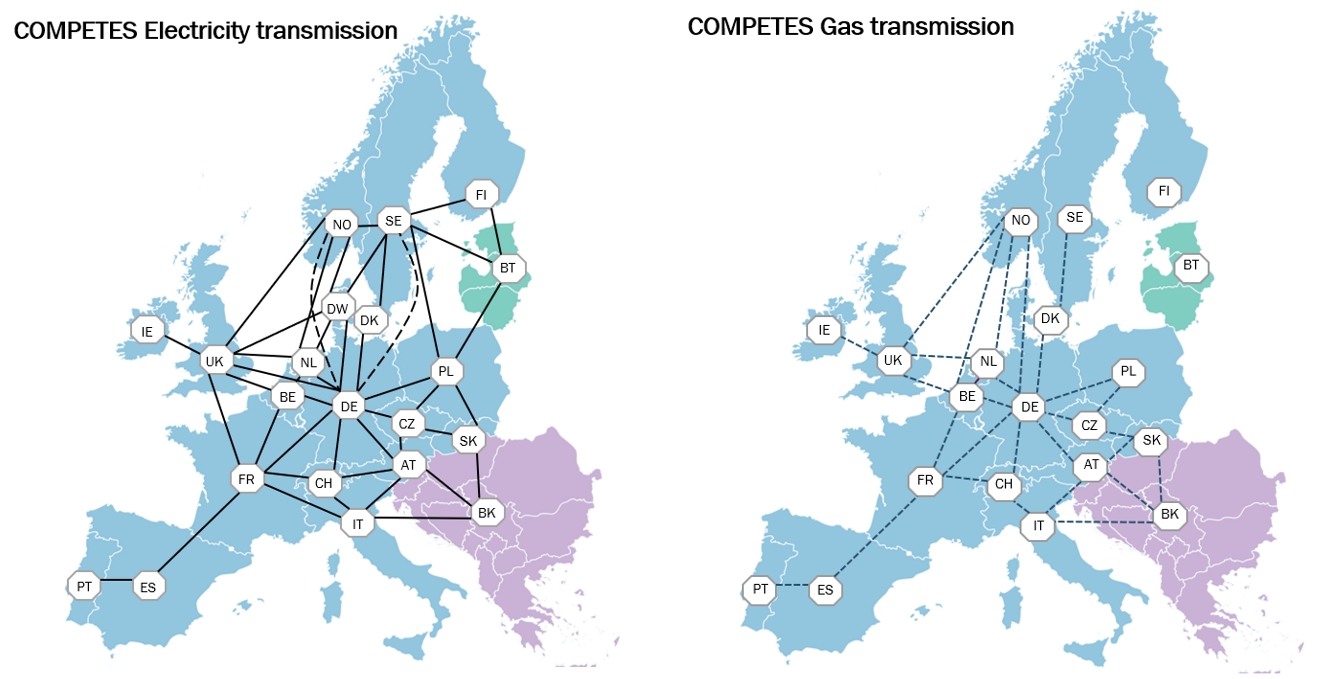}
    \caption{The geographical coverage of the COMPETES model}
    \label{fig:COMPETES_geogr_cov}
\end{figure}

The COMPETES model covers all EU Member States and some non-EU countries – i.e. Norway, Switzerland, the UK and the Balkan countries (grouped into a single Balkan region) – including a representation of the cross-border electricity transmission capacities interconnecting these European countries and regions; see Figure~\ref{fig:COMPETES_geogr_cov}. The model runs hourly, i.e., optimising the European power system over 8760 hours per year.

%%\FloatBarrier

Over the past two decades, COMPETES has been used for several assignments and studies on the Dutch and European electricity markets. In addition, it is used and regularly updated as part of the energy modelling framework for the annual Climate and Energy Outlook of the Netherlands; see, for instance, \citep{PBL2021a}.
For each scenario year, the primary inputs of COMPETES include parameters regarding the following features:
\begin{itemize}
    \item Electricity demand across all European countries/regions, including inelastic demand and additional demand due to further sectoral electrification of the energy system employing power-to-x technologies;
    \item Generation technologies, transmission interconnections, and flexibility options, including their techno-economic characteristics;
    \item Hourly profiles of various electricity demand categories and renewable energy technologies (notably hydro and variable renewable energy (VRE) sources such as solar and wind), including the full-load hours of these technologies;
    \item Expected future fuel and $CO_{2}$ prices;
    \item Policy targets/restrictions, such as meeting specific renewable energy or greenhouse gas (GHG) targets, or prohibiting the use of certain technologies, such as coal, nuclear, or Carbon Capture and Storage (CCS).
\end{itemize}
As indicated above, COMPETES includes a variety of flexibility options. More specifically, these options include:
\begin{itemize}
    \item Flexible generation: Conventional (gas, coal, nuclear) and renewable (curtailment of solar/wind);
    \item Cross-border electricity and $H_2$ trade;
    \item Demand response: Power-to-Mobility (P2M): electric vehicles (EVs), including grid-to-vehicle (G2V) and vehicle-to-grid (V2G); Power-to-Heat (P2H): industrial (hybrid) boilers and household (all-electric) heat pumps; Power-to-Gas (P2G), notably power-to-Hydrogen (P2H2);
    \item Storage: Pumped hydro (EU level), Compressed air (CAES/AA-CAES), Batteries (EVs, Li-ion, Lead-acid, Vanadium Redox), Underground storage of P2H2.
\end{itemize}

See \citet{sijm2017a} and \citet{ozdemir2020a} for a more detailed description of the COMPETES model. The explicit mathematical formulation of the optimisation model in COMPETES is available in \citet{ozdemir2019a}. It is important to note that the ramping constraints were not considered in this study for the sake of simplicity. However, the optimisation model does include demand response modelling based on the proposed formulation by \citet{moralesespana2022a}. In addition, a novel mathematical formulation to model the retrofit of the natural gas infrastructure to transport exclusively $H_2$ is included, as shown in Section \ref{sec:retrofit_formulation}. The enhanced version of the optimisation model allows a broader perspective on the possibilities of integrating $H_{2}$ in the 2050 reference scenario and its variants. The main scenario input parameters used for these scenarios are shown in the following section \ref{sec:refer_sc}.
\FloatBarrier

%---------------------------------------------------------------------------
%                Reference Scenario
%---------------------------------------------------------------------------
\subsection{Reference scenario}
\label{sec:refer_sc} 
In this section, we will discuss the supply and demand components for the 2050 European reference scenario and its variations. It is important to note that we specifically focus on modelling the electricity sector and other sectors which can potentially be electrified, such as transport, heat, and specifically $H_2$. However, the scenario does not consider the other potential uses of the natural gas sector in that year. This analysis assumes that the current natural gas network will be retrofitted to transport only $H_2$ or to meet any remaining natural gas demands.

\FloatBarrier
\subsubsection{Electricity demand}
\label{subsubsec:refer_sc_ener_dem_elec}
\begin{figure}[h]
    \centering
    \includegraphics[width=\FigureScale\textwidth]{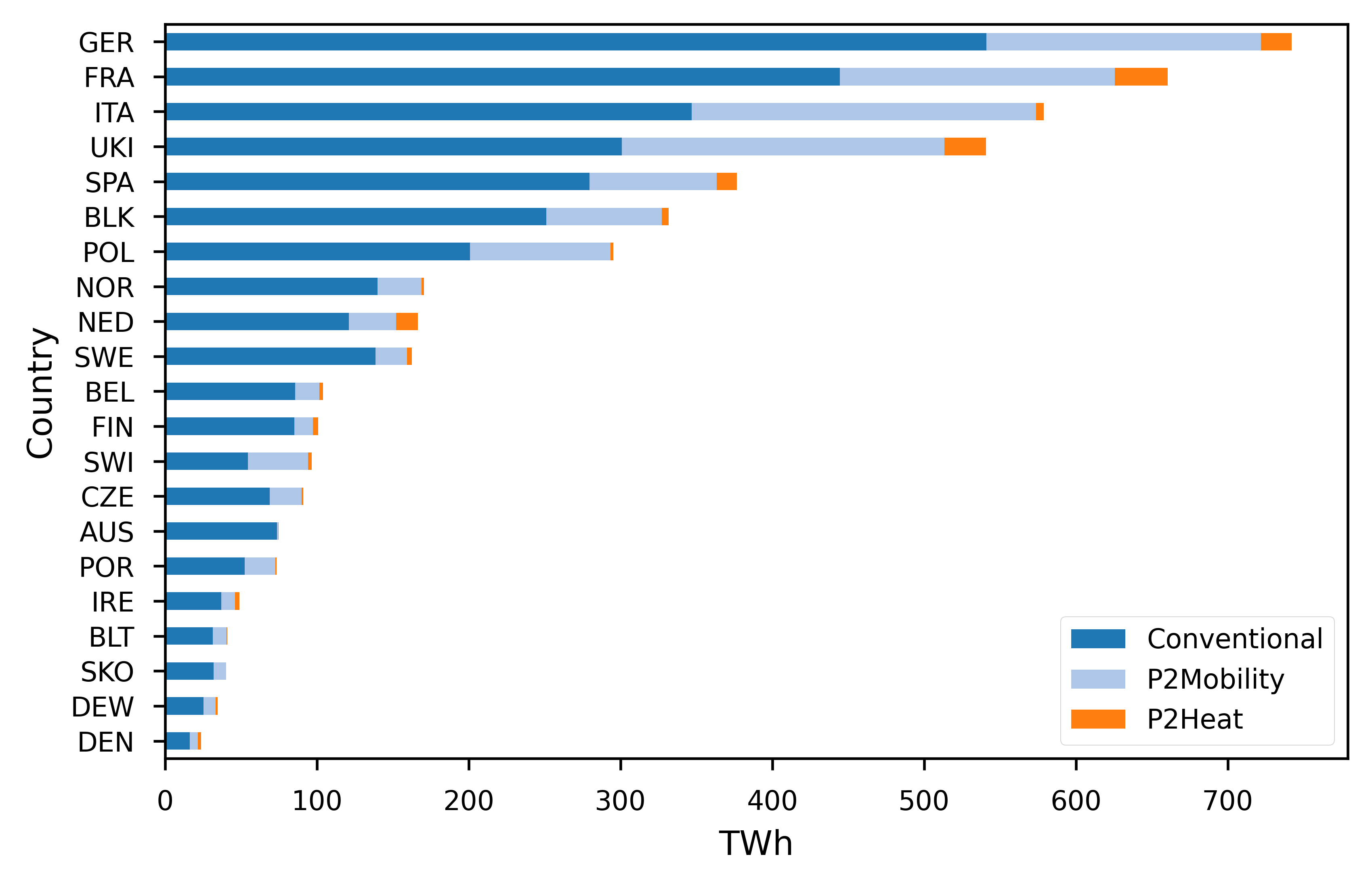}
    \caption{EU electricity demand per country}
    \label{fig:EU_elec_dem}
\end{figure}

Figure~\ref{fig:EU_elec_dem} provides an overview of the electricity demand parameters used in COMPETES for the EU countries in the reference scenario, referred to as 'R2050'. In this figure, the electricity demand is divided into three categories: Conventional, Power-to-Mobility, and Power-to-heat.The input parameters for the conventional demand category are incorporated into the model. Conversely, the model determines the Power-to-Mobility and Power-to-heat categories' final values endogenously. Below, there is additional information regarding each category.

\begin{itemize}
    \item Conventional electricity demand - For the R2050 scenario, the figures assume that the traditional demand for electricity consumption growth is offset more or less equally by the energy efficiency improvements. The hourly profile and demand per country are based on historical demand values \citep{entsoe2018a}.
    \item Power-to-Mobility – this demand for electric passenger vehicles (EVs) is assumed to be flexible to a certain extent. This demand includes both directions – i.e. grid-to-vehicle (G2V) and vehicle-to-grid (V2G). The projections on EV passenger vehicles are based on \citet{entsoe2018a} and \citet{berenschot2020a}. In addition, the EVs flexibility is endogenously considered in the model based on the modelling shown in \citep{moralesespana2022a}.
    \item Power-to-heat by households – This demand comes from electric heat pumps; similar to EVs, this demand is assumed to be flexible and determined endogenously using the formulation proposed in \citep{moralesespana2022a}. A set of constraints limits the flexibility of the heat pump.  The projections on household heat pumps are based on \citet{entsoe2018a} and \citet{berenschot2020a}.
\end{itemize}

\FloatBarrier
\subsubsection{Hydrogen demand}
\label{subsubsec:refer_sc_ener_dem_H2}
\begin{figure}[h]
    \centering
    \includegraphics[width=\FigureScale\textwidth]{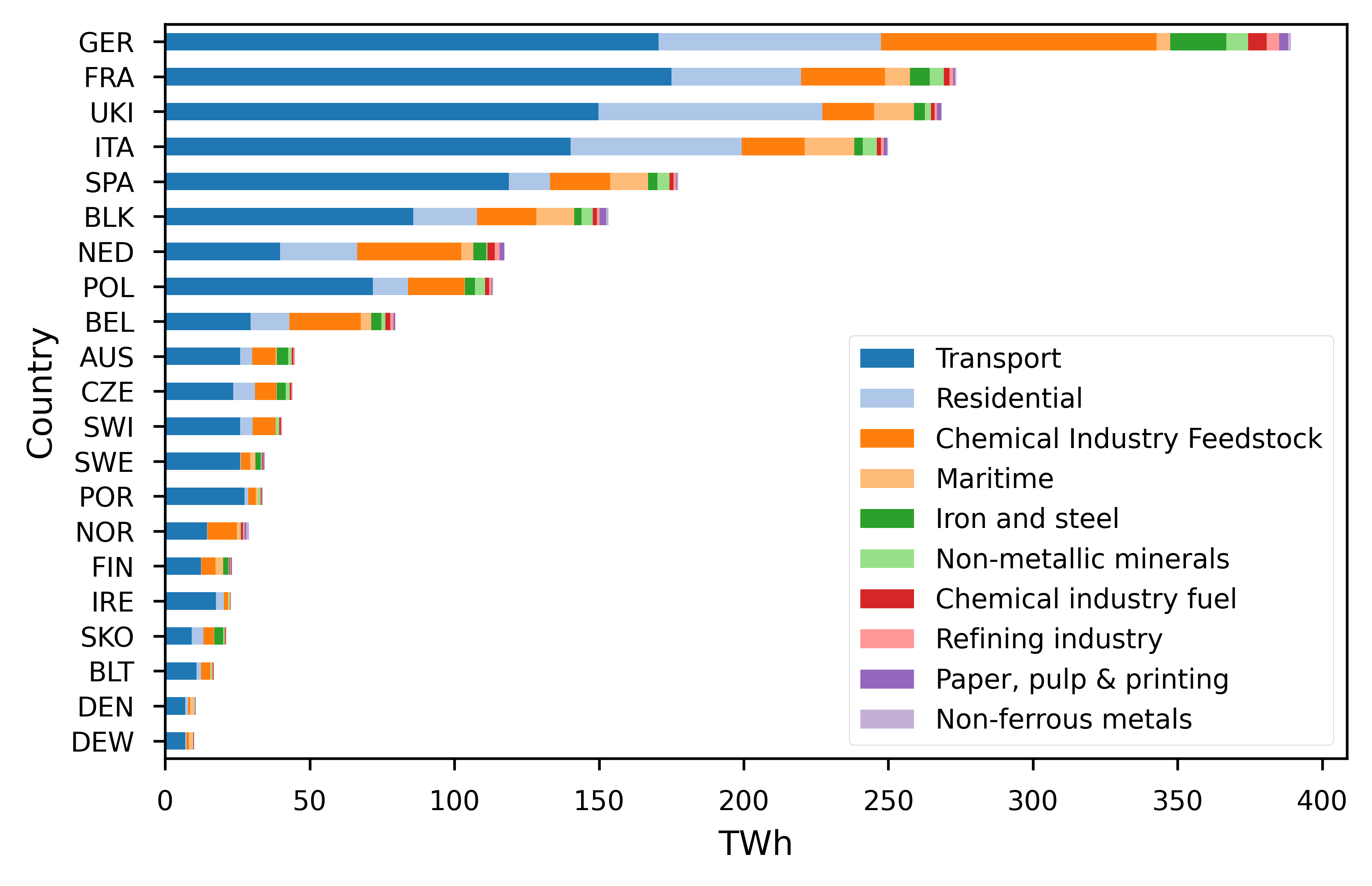}
    \caption{EU $H_2$ demand per country and sector}
    \label{fig:EU_H2_dem}
\end{figure}

The $H_2$ demand is based on one of the eight scenarios from the European Commission's long-term strategy to reduce greenhouse gas.
\textit{The Commission's analysis is based on the PRIMES, GAINS, GLOBIOM model suite and explores eight economy-wide scenarios to achieve different levels of ambition for 2050, covering the potential range of reduction needed in the EU to contribute to the Paris Agreement's temperature objectives of between the well below 2°C and to pursue efforts to limit to 1.5°C temperature change} \citep{comission2018a}.
The selected scenario, '1.5TECH' described in \citet{comission2018a}, focuses on technical solutions to achieve net-zero GHG emissions. It increases CCS and uses E-gases and E-fuels based on air-capture or biogenic $CO_2$ to further reduce emissions. The scenario also applies negative emission technologies via biomass coupled with CCS and the storage of biogenic $CO_{2}$. Figure~\ref{fig:EU_H2_dem} shows the $H_2$ demand per EU country and for different activities within the industrial sector.

\FloatBarrier
\subsubsection{Energy supply: sources and technologies}
\label{subsec:refer_sc_ener_supply}
\begin{figure}[h]
    \centering
    \includegraphics[width=\FigureScale\textwidth]{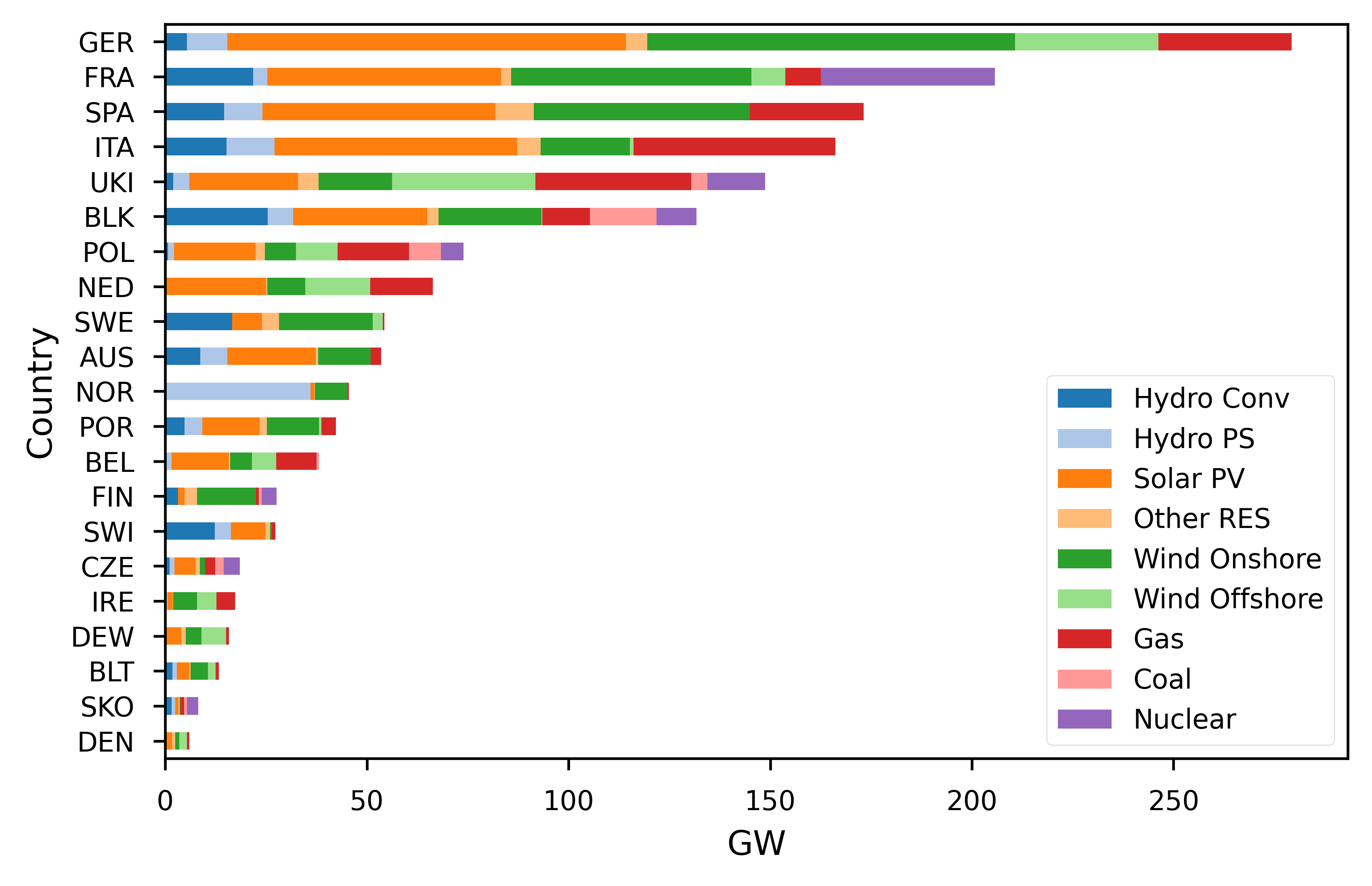}
    \caption{Initial installed capacities per country}
    \label{fig:EU_init_installed_cap}
\end{figure}

COMPETES uses its investment module to meet the demand of the different energy vectors, i.e. electricity and $H_2$, in a cost-optimal way. COMPETES includes a variety of primary energy sources and technologies and energy conversion technologies. These technologies are described in Section~\ref{sec:opt_model}. 
Figure~\ref{fig:EU_init_installed_cap} shows the initial electricity generation capacities in the COMPETES model. These initial capacities serve as an input for the model, which will help to determine the new required capacity to meet the 2050 electricity and $H_2$ demand. These initial input values are based on the National Trends scenario \citep{entsoe2020a}.

\FloatBarrier
\subsubsection{Hydrogen generation}
\label{subsubsec:refer_sc_H2_tech_gen}
Similar to electricity generation, initial $H_2$ generation values are defined exogenously in the model. The $H_2$ generation technologies in COMPETES are steam-methane reforming (SMR), SMR with a 54\% CCS rate (SMR CCS 54), SMR with an 89\% CCS rate (SMR CCS 89), and electrolysers. Table \ref{tab:electrolyser_smr_technical_data} displays the techno-economic parameters that were taken into account for the $H_2$ generation technologies. Among the available options for electrolysers, alkaline electrolysis (AE) and proton-exchange membrane (PEM) electrolysis are the leading technologies. PEM has advantages like a faster ramping rate and is projected to become more cost-effective in the future, as explained in \citep{BOHM2020}. As stated in Section \ref{sec:opt_model}, the COMPETES model does not take into account ramping constraints for the sake of simplicity. As a result, the electrolyser technology listed in Table \ref{tab:electrolyser_smr_technical_data} can serve as a representation of both types of technologies based on the CAPEX and Fixed O\&M assumptions, as the difference in efficiency is minimal in the year 2050 \citep{BOHM2020}.

\begin{table}[h]
    \centering
    \begin{tabular}{ccccc}
        \hline
        Technology&SMR&SMR  &SMR  &Electrolyser  \\
                  &   &CCS 54&CCS 89&              \\        
        \hline
        Source (X)                &Gas  &Gas  &Gas  &Power \\
        2use (2Y)                 &$H_2$&$H_2$&$H_2$&$H_2$ \\
        Capex [\euro{}/kW]        &744  &881  &1330 &600   \\
        Fixed O\&M [\euro{}/kW/yr]&27   &44   &62   &20    \\
        LifeTime [Years]          &25   &25   &25   &30    \\
        Efficiency [p.u.]         &0.76 &0.74 &0.69 &0.68  \\
        Emissions [kg/MWh]        &229  &105  &26   &0     \\
        CCS [kg/MWh]              &-    &124  &204  &-     \\
        \hline
    \end{tabular}
    \caption{Techno-economic input parameters for $H_{2}$ technologies}
    \label{tab:electrolyser_smr_technical_data}
\end{table}
%%\FloatBarrier

Table \ref{tab:electrolyser_smr_cap} presents the initial $H_{2}$ generation output capacities assumed in this study. These are based on \citet{maisonnier2007a} and \citet{fticonsulting2020a}. It is worth noting that the latest published ambitions of the European countries are higher than these values.

\begin{table}[h]
    \centering
    \begin{tabular}{ccc}
        \hline
        Country & Electrolyser [MW] & SMR [MW] \\
        \hline
        GER & 1000 & 1900 \\
        FRA & 6500 & 530 \\
        UKI & 5000 & 144 \\
        SPA & 4000 & 554 \\
        NED & 3000 & 1144 \\
        POR & 2000 & 25 \\
        AUS & 0 & 90 \\
        BEL & 0 & 783 \\
        SWI & 0 & 28 \\
        CZE & 0 & 141 \\
        DEN & 0 & 29.5 \\
        DEW & 0 & 29.5 \\
        FIN & 0 & 413 \\
        BLK & 0 & 523 \\
        IRE & 0 & 0 \\
        ITA & 0 & 411 \\
        NOR & 0 & 0 \\
        SWE & 0 & 0 \\
        SKO & 0 & 115 \\
        POL & 0 & 0 \\
        BLT & 0 & 221 \\
        \hline
    \end{tabular}
    \caption{Initial $H_2$ output capacities}
    \label{tab:electrolyser_smr_cap}
\end{table}

\FloatBarrier
\subsubsection{Hydrogen transport}
\label{subsubsec:refer_sc_H2_tech_transp}
COMPETES recently introduced the possibility of investing in $H_2$ pipelines for long-distance, high-volume $H_2$ transport between European countries, considering a transport model. Figure~\ref{fig:COMPETES_geogr_cov} shows the existing natural gas transport links, which can be repurposed to deliver $H_2$ instead of natural gas. Also, the model can invest in new $H_2$ pipelines. The following decisions can be made endogenously by the model: invest in new $H_2$ pipelines, retrofit to 60\% of initial gas capacity, retrofit to 80\% of initial gas capacity (more expensive than 60\% due to extra compression needed). Moreover, the modelling assumptions include:
    \begin{itemize}
        \item Retrofit and investment decisions are only possible where there is an existing pipeline.
        \item Transport capacities are assumed to be bidirectional, e.g. assume that gas trading capacity between Germany and the Netherlands is different depending on the trade direction. In these cases, we consider the highest capacity for both directions.
        \item Neither losses in transport nor variable costs for compression are considered in this study. for this study. Future research conducting a more in-depth analysis of the impacts of pressures and compression losses may offer valuable and informative perspectives on this subject.
        \item In the $H_2$ sector, decision variables are used every six hours, while the power system variables are decided upon hourly. This temporal difference in the variables' definition allows the model to account for the storage and time-shifting of $H_2$ through pipelines in a simplified way.        
        \item Pipelines can be retrofitted partially since the model is only able to make continuous investment decisions.
    \end{itemize}

Furthermore, it is assumed that SMR can supply a maximum of 50\% of the $H_{2}$ demand of a country; this is enforced with a constraint in COMPETES. This is based on the RED II revision proposal \citep{comission2021a}, which sets a binding 50\%\footnote{The European Commission has recently increased its target to 70\% by 2035, according to the REPowerEU plan.} target for renewable fuels of non-biological origin used as feedstock or energy carriers.

The natural gas transmission network consists of multiple parallel pipelines. Therefore, the results for $H_{2}$ transport in this paper can be interpreted as the number of pipelines that require retrofitting for $H_{2}$ transportation instead of natural gas.

\FloatBarrier
\subsubsection{Fuel and $CO_{2}$ prices}
\label{subsec:refer_sc_CO2_prices}
Table \ref{tab:fuel_CO2_prices} shows the fuel prices taken from references \citep{berenschot2020a} and \citep{PBL2021a}. The $CO_{2}$ price was considered as 250~\euro{}/ton. It is important to highlight that if natural gas prices increase beyond what is listed in this table, it could encourage the electricity sector to produce more $H_2$ using renewable energy sources.
\begin{table}[h]
    \centering
    \begin{tabular}{cc}
        \hline
        Fuel & Price 2050 [\euro2015/GJ] \\
        \hline
        Oil & 10.63 \\
        Biomass & 9.00 \\
        Natural Gas & 7.54 \\
        Coke Oven Gas & 7.54 \\
        Coal & 2.25 \\
        Lignite & 1.10 \\
        Nuclear & 0.78 \\
        \hline
    \end{tabular}
    \caption{Fuel prices}
    \label{tab:fuel_CO2_prices}
\end{table}

\FloatBarrier
\subsubsection{Carbon capture and storage (CCS)}
\label{subsec:refer_sc_CCS}
COMPETES endogenously optimises the investments in electricity and $H_2$ generation units with carbon dioxide capture and storage (CCS), such as: Biomass plants with CCS, gas CCGT plants with CCS, coal-fired plants with CCS, SMR CCS 54, and SMR CCS 89.

Importantly, $CO_{2}$ geological storage is currently prohibited in some countries. This study uses current national legislations and regulations to determine whether the model can invest in the aforementioned technologies. Based on the EU Directive 2009/31/EC on the geological storage of $CO_{2}$ \citep{CO2GeoNet2021a}: Germany, Austria, Estonia, Latvia, Lithuania, Denmark, Finland, and Ireland do not allow $CO_{2}$ geological storage.

\FloatBarrier
%---------------------------------------------------------------------------
%           Impact of H2 electrification
%---------------------------------------------------------------------------
\section{Results of the impact of $H_{2}$ electrification}
\label{sec:H2electr} 

\subsection{Reference scenario and the impact of $H_{2}$ electrification}
\label{subsec:H2electr_RefScn}
This section provides the main results of the reference scenario described in Section~\ref{sec:refer_sc} and the impact of $H_{2}$ electrification. To measure the effect of $H_{2}$ electrification, the reference scenario R2050 is compared with scenario NoP2H2, where electrolysis is not allowed, i.e., all $H_{2}$ demand must be supplied via SMR (except for the initial electrolysers installed capacity).

\begin{figure}[h]
    \centering
    \includegraphics[width=\FigureScale\textwidth]{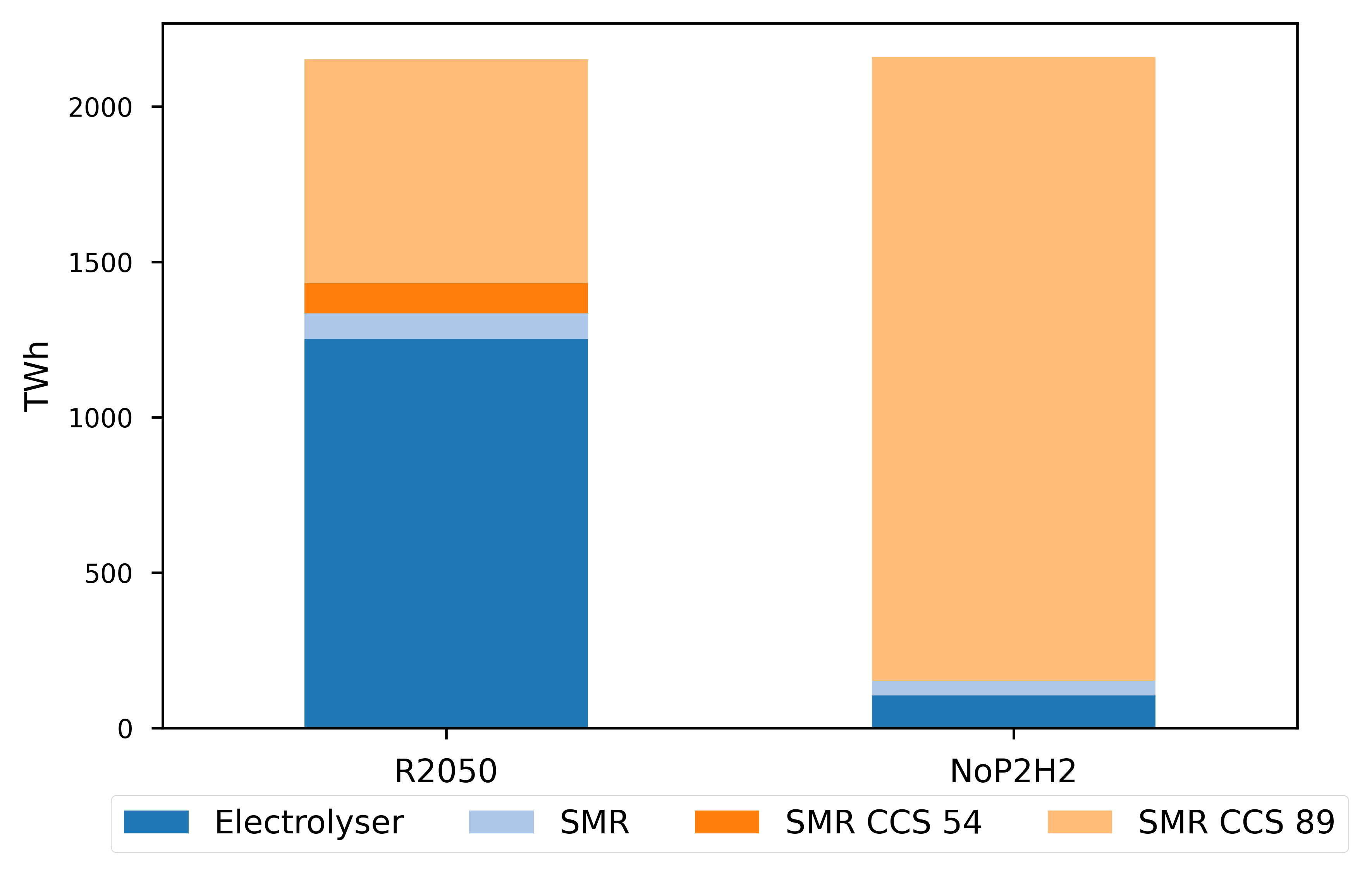}
    \caption{EU $H_2$ supply in R2050 and NoP2H2}
    \label{fig:EU_H2_supply}
\end{figure}

Figure~\ref{fig:EU_H2_supply} shows the $H_2$ supply sources for both scenarios, where the reference scenario R2050 electrifies 58\% of the total $H_{2}$ demand, and the remaining 42\% is supplied via SMR. The SMR technology that dominates is SMR with CCS with an 89\% $CO_{2}$ capture rate (SMR CCS 89) due to the $CO_{2}$ price of 250~\euro{}/ton. Notice that 5\% of the $H_{2}$ demand in NoP2H2 is provided by electrolysis, which results from the initially installed capacities of electrolysers. Similarly, the model uses existing SMR without CCS facilities to supply 4\% and 2\% of the $H_{2}$ demand in R2050 and NoP2H2, respectively. 

\begin{figure}[h]
    \centering
    \includegraphics[width=\FigureScale\textwidth]{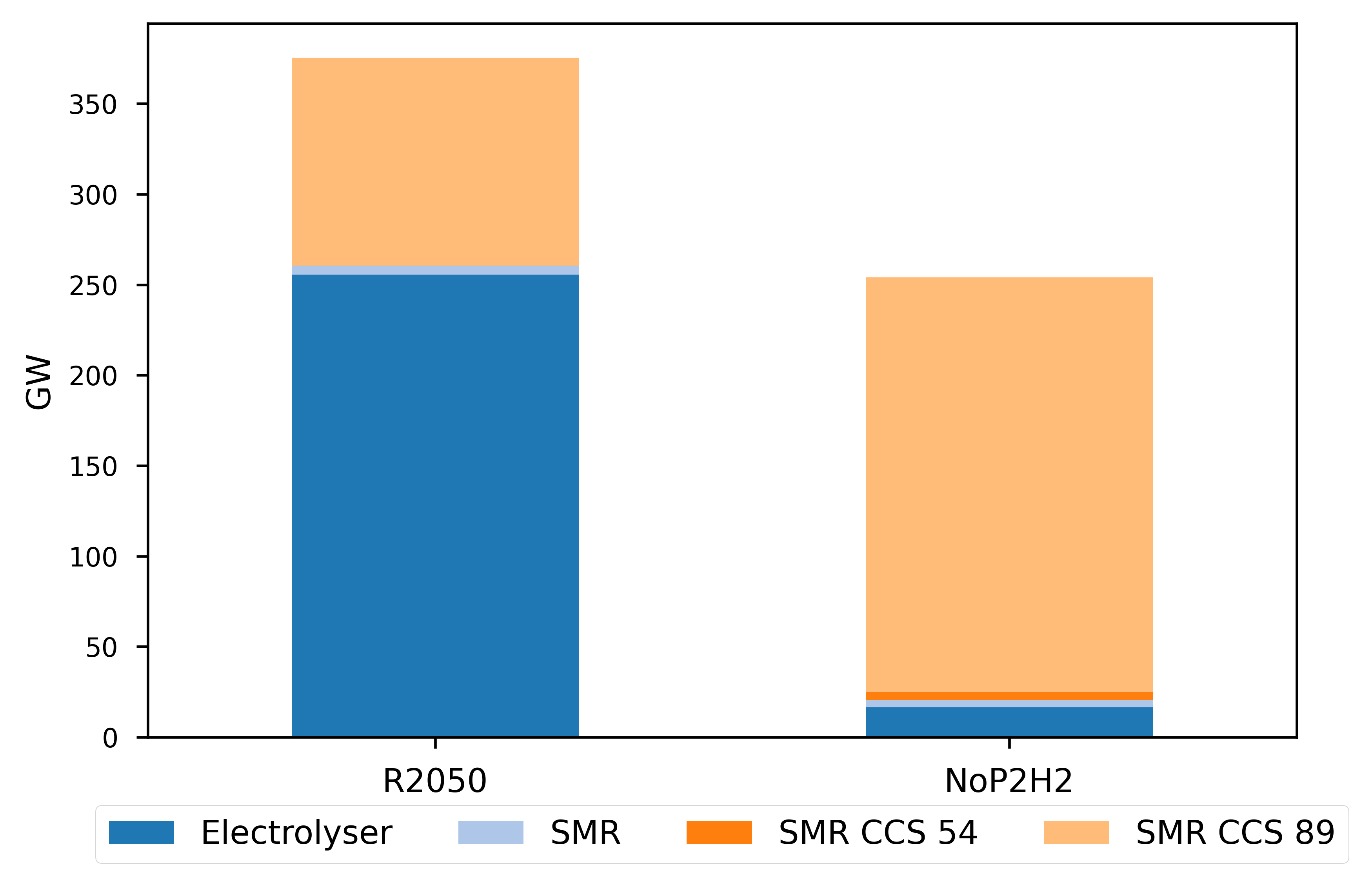}
    \caption{EU $H_{2}$ output capacity in R2050 and NoP2H2}
    \label{fig:EU_H2_capacity}
\end{figure}

Figure~\ref{fig:EU_H2_capacity} shows the total installed capacities for $H_{2}$ supply. Notice that the total $H_{2}$ installed capacity of R2050 increases 48\% compared to NoP2H2. This capacity oversize in the R2050 is caused by the investment in electrolysers, which requires the installation of approximately 250 GW capacity in the EU for the year 2050. This result suggests that it is more economically efficient to invest in a larger capacity to produce more $H_{2}$ during periods with low electricity prices (VRE-dominated production) and store it for later use, thus avoiding producing $H_{2}$ during high electricity prices. Electrolysers present 4900 Full Load Hours (FLH) in R2050, whereas SMR operates 7500 FLH in R2050 and 8650 FLH in NoP2H2.

\FloatBarrier

Table~\ref{tab:countries_H2_supply} shows how the $H_{2}$ is supplied in different countries via SMR, electrolysis, or imports. One can observe three different types of countries:
\begin{enumerate}
    \item Countries with very high VRE production, like France and Spain, which shift most of their $H_{2}$ production from SMR in NoP2H2 to electrolysis in R2050, and even become net exporters.
    \item countries that still find SMR as the most economical way to produce $H_{2}$ (e.g., because of low VRE potentials), which use their maximum allowed SMR production, 50\% of the internal demand, and supply the remaining $H_{2}$ demand via electrolysis and imports, even moving from a net export position in NoP2H2 to a net import position; this is the case for countries like the Netherlands and the UK.
    \item Countries that do not allow carbon storage, which supply $H_{2}$ mainly through imports and electrolysis, e.g., Germany and Austria.
\end{enumerate}

\begin{table}[h]
    \centering
    \begin{tabular}{cccccc}
        \hline
        Country   & Scenario   &   Electrolyser &   SMR       &   SMR       &   Net        \\
                  &            &                &   CCS 54    &   CCS 89    &   import     \\        
        \hline
        NED       & R2050      &             23 &           2 &          57 &            0 \\
        NED       & NoP2H2     &              9 &           1 &         417 &         -310 \\
        GER       & R2050      &             75 &           5 &           0 &          193 \\
        GER       & NoP2H2     &             30 &           2 &           0 &          363 \\
        AUS       & R2050      &             29 &           0 &           0 &           21 \\
        AUS       & NoP2H2     &              0 &           0 &           0 &           45 \\
        FRA       & R2050      &            348 &           1 &         136 &          -83 \\
        FRA       & NoP2H2     &             25 &           1 &         251 &           -3 \\
        UKI       & R2050      &             69 &           0 &         134 &            1 \\
        UKI       & NoP2H2     &             17 &           0 &         273 &          -22 \\
        SPA       & R2050      &            160 &           1 &          88 &          -27 \\
        SPA       & NoP2H2     &             16 &           1 &         159 &            2 \\
        \hline
    \end{tabular}
    \caption{$H_2$ supply in R2050 and NoP2H2 [TWh]}
    \label{tab:countries_H2_supply}
\end{table}

\begin{figure}[h]
    \centering
    \includegraphics[width=\FigureScale\textwidth]{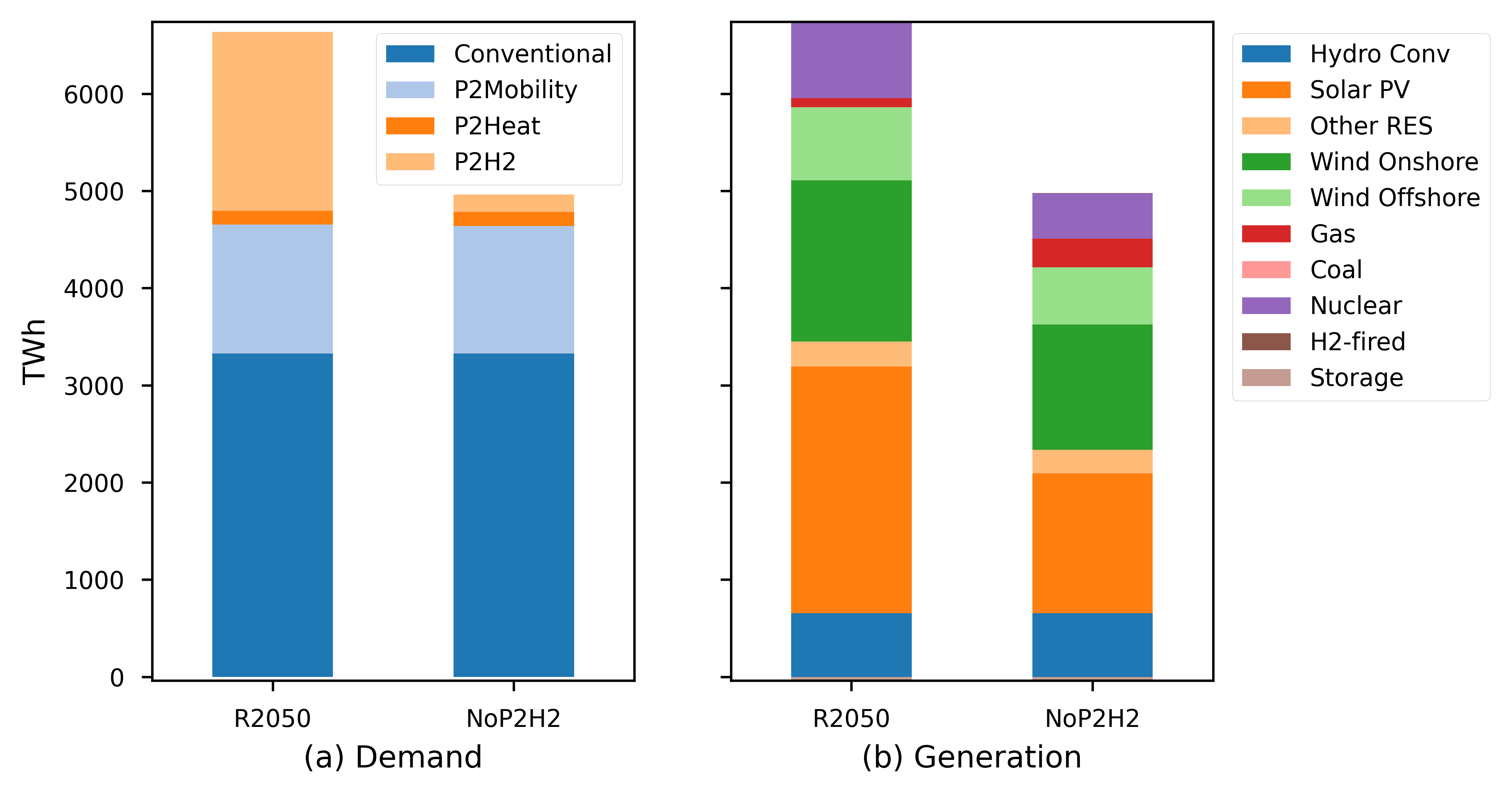}
    \caption{EU electricity demand and supply in R2050 and NoP2H2}
    \label{fig:EU_power_demand_supply}
\end{figure}

\begin{figure}[h]
    \centering
    \includegraphics[width=\FigureScale\textwidth]{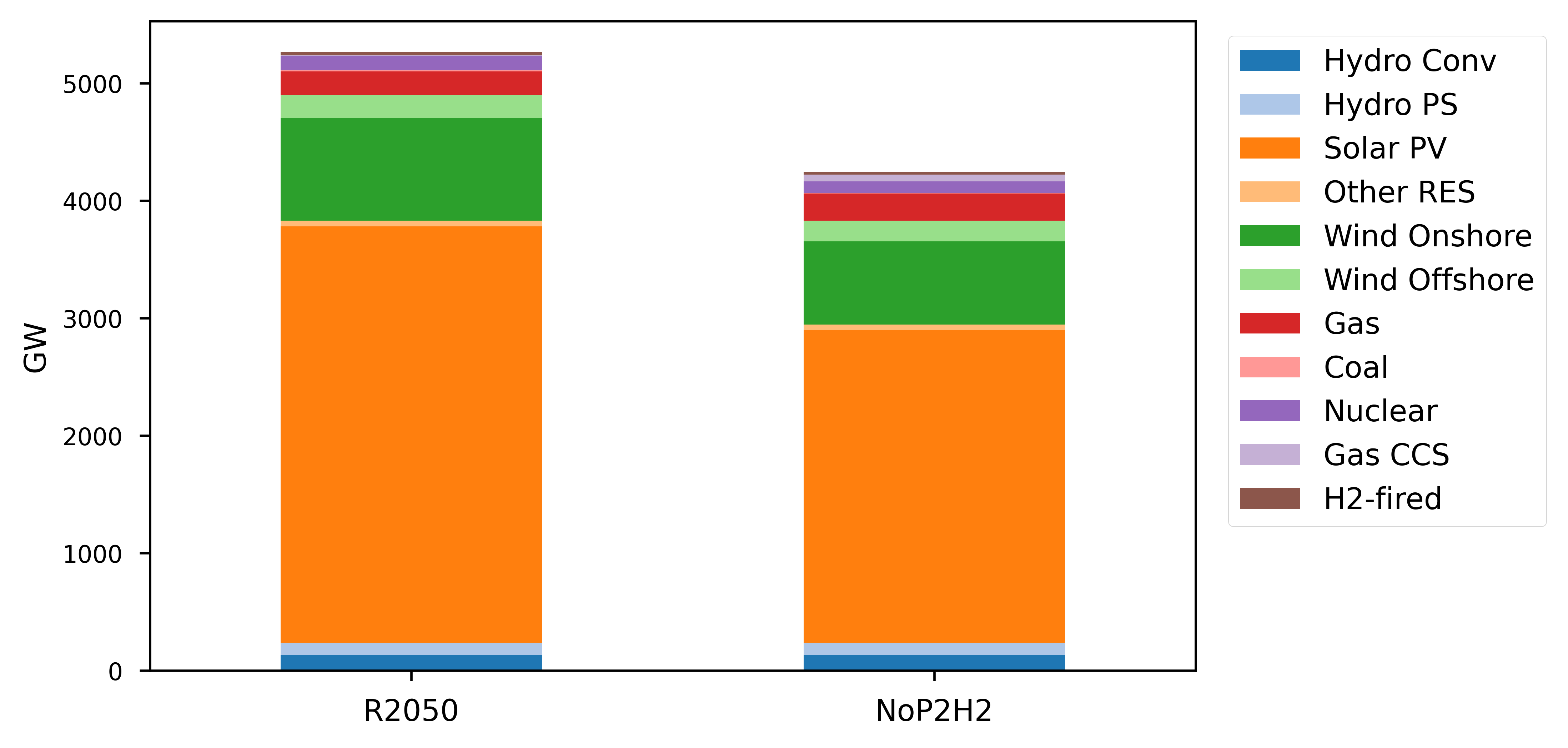}
    \caption{EU Generation capacity in R2050 and NoP2H2}
    \label{fig:EU_generation_capacity}
\end{figure}

The 58\% level of $H_{2}$ electrification requires almost 1841~TWh\footnote{The EU total electricity demand in 2021 was 2865~TWh.} of extra electricity demand, accounting for nearly 28\% of the total electricity demand in 2050. See Figure~\ref{fig:EU_power_demand_supply}, which shows the electrical energy mix for both scenarios R2050 and NoP2H2. To meet the new $H_{2}$ electrical demand, the system uses and invests mainly in more solar PV, wind offshore, and wind onshore; see Figure~\ref{fig:EU_generation_capacity}. The total VRE production changes from 3316 TWh in NoP2H2 to 4952 TWh in R2050. This VRE production increase is equivalent to 90\% of the extra $H_{2}$ electricity demand. Nuclear production also increases from 463 TWh to 789 TWh, equaling 18\% of the additional $H_{2}$ electrical demand. The production of non-polluting technologies, VRE and nuclear, coveres the extra $H_{2}$ demand and even replaces part of the gas production, which decreases from 295 TWh in NoP2H2 to 93 TWh in R2050. This reduction of almost 70\% in gas production appears due to the flexibility offered by the electrolysers through shifting in time (storage) and shifting source (to SMR), where electrolysis increases when electricity prices are low (e.g., due to VRE abundance). This allows the new VRE and nuclear investments to be better used, and when electricity prices are high (e.g., due to production of gas-fired power plants) electrolysis may not be viable. However, the extra investment in VRE and nuclear is still present, thus replacing gas sources.

Figure~\ref{fig:EU_generation_capacity} shows that there is significant increase of VRE capacity to cover the new $H_{2}$ electrification in R2050. Moreover, the flexibility of electrolysers allows for relying less on peak units. For instance, the installed capacity of gas-CCS power plants reduces from 58 GW in NoP2H2 to 6.50 GW in R2050.

\FloatBarrier

\begin{figure}[h]
    \centering
    \includegraphics[width=\FigureScale\textwidth]{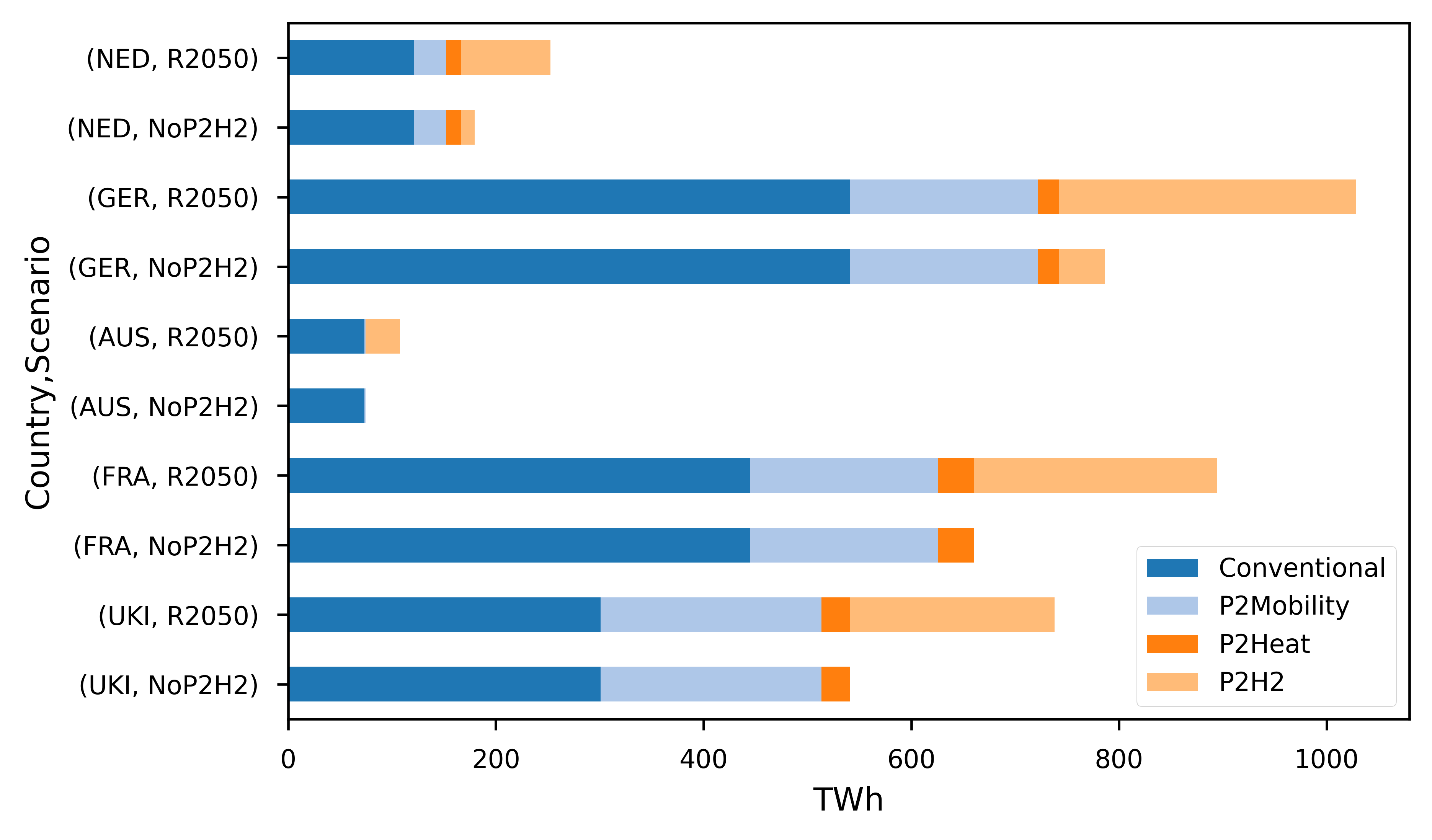}
    \caption{NED, GER, AUS, FRA, UKI and SPA electricity} demand in R2050 and NoP2H2
    \label{fig:countries_power_demand}
\end{figure}

\begin{figure}[h]
    \centering
    \includegraphics[width=\FigureScale\textwidth]{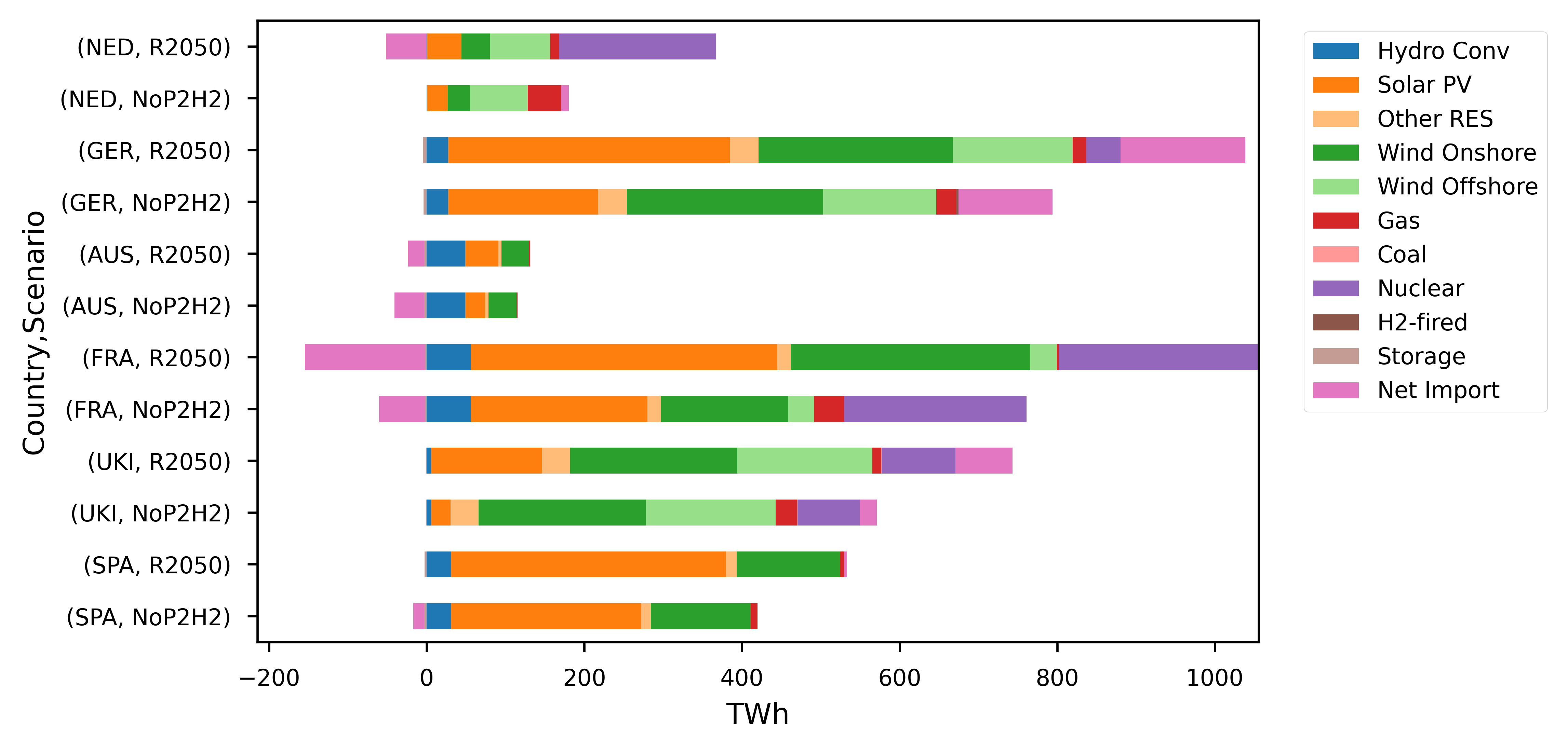}
    \caption{NED, GER, AUS, FRA, UKI and SPA electricity} supply in R2050 and NoP2H2
    \label{fig:countries_power_supply}
\end{figure}

Note that although some countries do not allow carbon storage, they still import $H_{2}$, which is mainly produced via SMR with CCS in NoP2H2. Some of these countries even use this SMR-generated $H_2$ to generate electricity--as is the case for Germany, producing 3 TWh via $H_2$-2-power--as shown in Figures \ref{fig:countries_power_demand} and \ref{fig:countries_power_supply}. Furthermore, it is possible that even in the case of R2050, these countries still import $H_{2}$ that is produced with SMR with CCS; once the $H_{2}$ is produced and injected into a pipeline, it is not possible to know if some given molecules of $H_{2}$ were produced via electrolysis or SMR. Figures~\ref{fig:countries_power_demand} and \ref{fig:countries_power_supply} also show how the demand and supply of electricity are divided within different EU countries. As expected, the countries with the highest electrolysis also have the highest VRE production. For example, France, Spain, and Germany, with electrolysis demand increases of 198 TWh, 127 TWh, and 249 TWh, saw increases in VRE production of 398 TWh, 138 TWh, and 170 TWh, respectively, in R2050 compared to NoP2H2.

\FloatBarrier

\begin{figure}[h]
    \centering
    \includegraphics[scale=0.4]{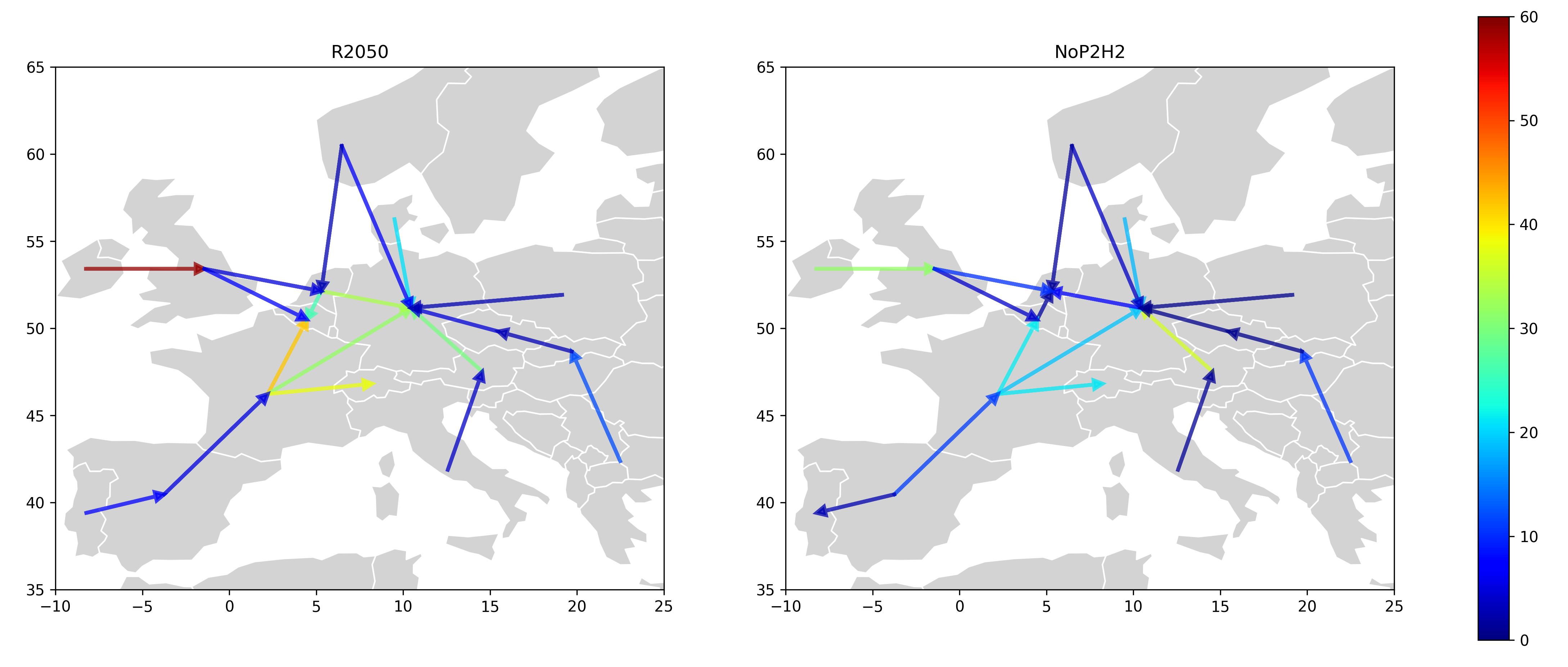}
    \caption{Electricity net trade patterns [TWh]}
    \label{fig:net_elec_trade}
\end{figure}

Table~\ref{tab:EU_elec_transm_invest} shows the new electricity and $H_{2}$ transmission investments for R2050 and NoP2H2. Interestingly, R2050 requires almost a 45\% lower electricity transmission capacity, even though its electricity demand is 35\% higher than NoP2H2. This is caused by the flexibility offered by the electrolysers, where the $H_{2}$ electrification helps the system rely less on other sources of flexibility, such as electricity trade and peak units. Figure~\ref{fig:net_elec_trade} shows this general pattern of higher electricity trade in R2050; countries presenting VRE abundance trade more, as is the case of Portugal and Spain exporting more to France, where Portugal changes from being a net importer in NoP2H2 to a net exporter in R2050. Another country that becomes a net exporter is the Netherlands, not only because of its extra VRE production, but also its extra nuclear production; see Figure~\ref{fig:countries_power_supply}.

\begin{table}[ht]
    \centering
    \begin{tabular}{ccc}
        \hline
        [GW] & E-transmission & $H_{2}$ Transmission \\
        \hline
        R2050  & 28 & 90 \\
        NoP2H2 & 52 & 56 \\
        \hline
    \end{tabular}
    \caption{EU electricity and $H_{2}$ transmission investments in R2050 and NoP2H2}
    \label{tab:EU_elec_transm_invest}
\end{table}

\FloatBarrier

\begin{figure}[h]
    \centering
    \includegraphics[scale=0.4]{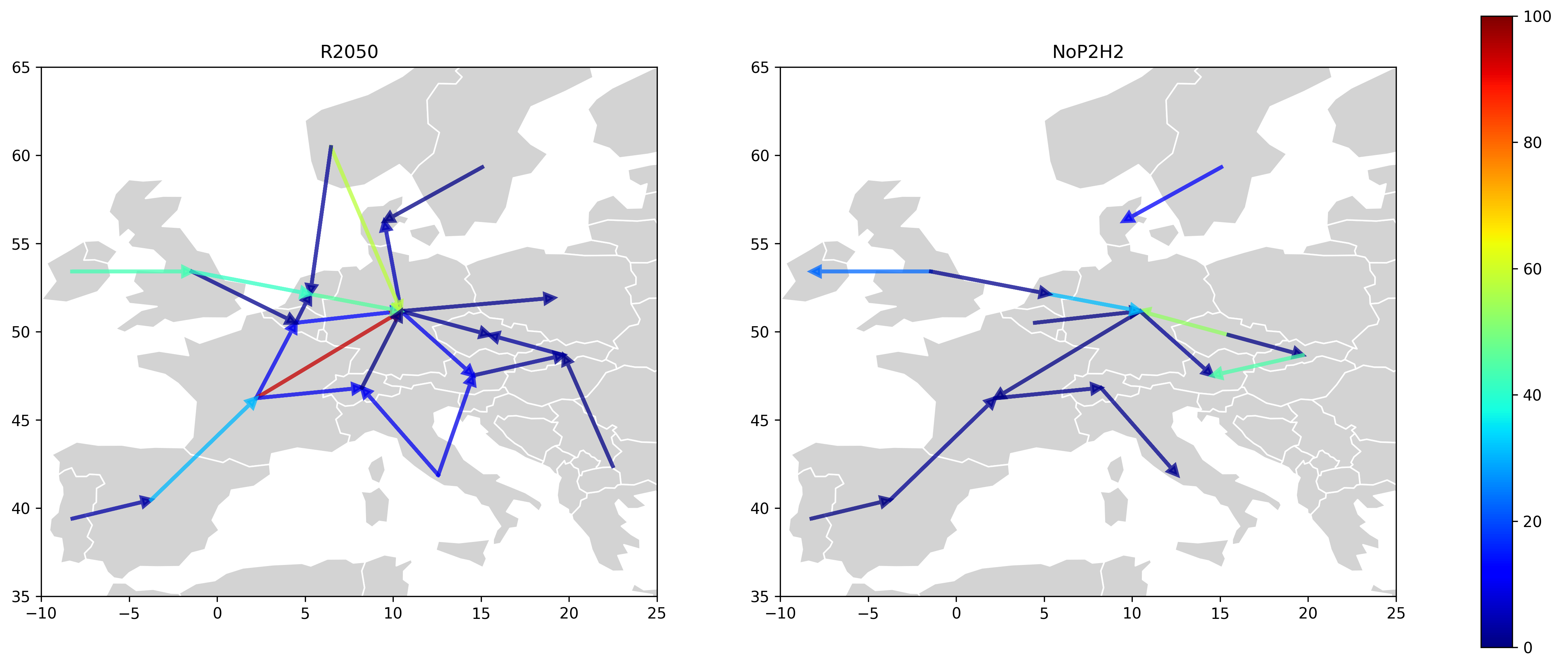}
    \caption{$H_{2}$ net trade patterns [TWh]}
    \label{fig:net_H2_trade}
\end{figure}

The $H_{2}$ transmission capacity is 1.6x higher in R2050 (see Table~\ref{tab:EU_elec_transm_invest}), resulting in  different trade patterns compared to NoP2H2; see Figure~\ref{fig:net_H2_trade}. This is a natural consequence of countries with higher VRE investments and generation, producing (and exporting) $H_2$ from electrolysis. For example, France, Spain, and Norway (with high solar and wind) export more $H_2$ to Germany (which has high import demand due to CCS facilities not being allowed). Ireland also supplies the UK with higher offshore wind. The Netherlands shows a less significant trade volume change with Germany due to the 50\% SMR production restriction policy.

\begin{figure}[h]
    \centering
    \includegraphics[scale=0.5]{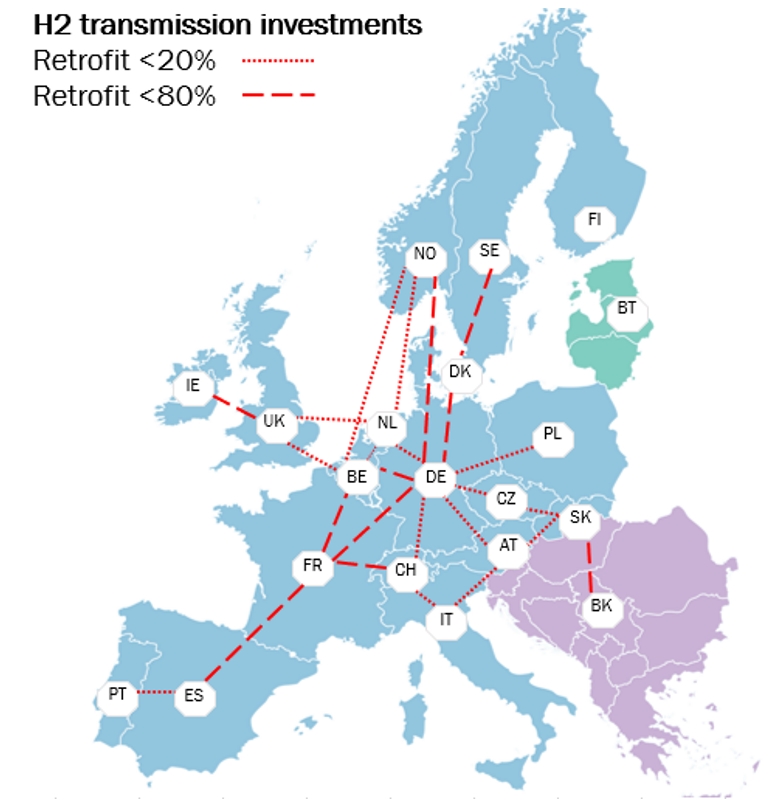}
    \caption{EU $H_{2}$ transmission investments map in R2050}
    \label{fig:EU_H2_trans_invest}
\end{figure}

The current gas infrastructure already offers enough potential to accommodate the need for $H_{2}$ trade. Figure~\ref{fig:EU_H2_trans_invest} shows the different $H_{2}$ transmission investments between countries where no new pipelines were built and 11\% of the existing gas infrastructure's total capacity was retrofitted for $H_{2}$ transportation. Although we did not model gas (methane) transportation, the remaining transport capacity (89\%) is enough to accommodate 2050 requirements, which will be lower than current requirements. The NoP2H2 scenario only needed to retrofit 7\% of the gas infrastructure's total capacity due to the lower need for $H_{2}$ trade. 

\FloatBarrier

\begin{figure}[h]
    \centering
    \includegraphics[width=\FigureScale\textwidth]{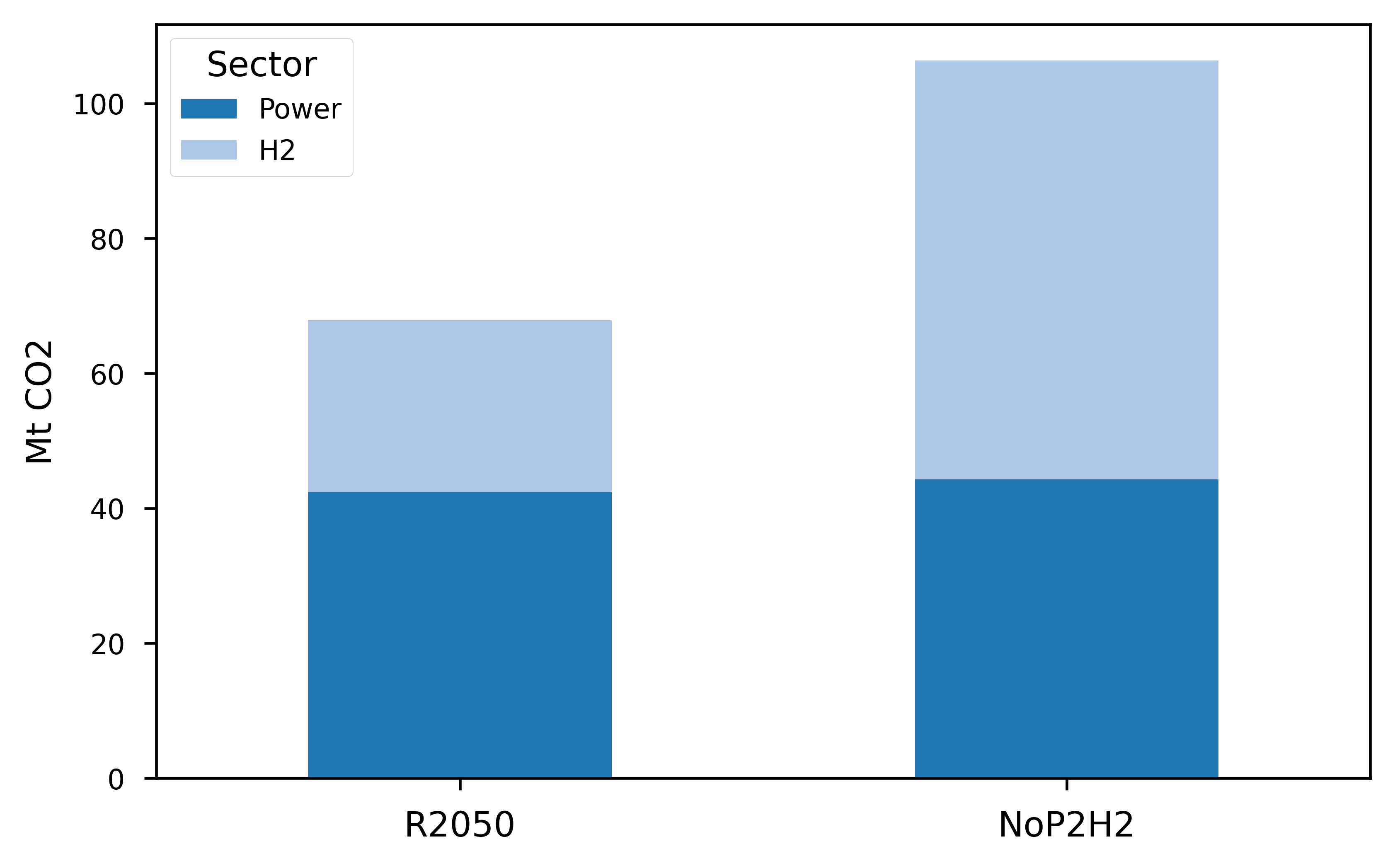}
    \caption{EU $CO_{2}$ emissions for the electricity} and $H_{2}$ sector in R2050 and NoP2H2
    \label{fig:EU_CO2_emissions}
\end{figure}

The total $CO_{2}$ emissions of the system decreased from 106 Mton in NoP2H2 to 68 Mton in R2050. This 35\% emissions reduction results from shifting 58\% of the $H_{2}$ production from SMR, in NoP2H2, to electrolysis, in R2050, where mainly non-pollutive (VRE and nuclear) technologies supply the electricity for electrolysis. Even though the power system has additional demand, its $CO_{2}$ emissions are lower in R2050 since non-polluting technologies are also replacing part of the gas-fired technologies that were present in NoP2H2. As shown in Figure~\ref{fig:EU_CO2_emissions}, the $CO_{2}$ emissions in the electricity sector decreased from 44~Mton in NoP2H2, to 42~Mton in R2050. In short, electrifying part of the $H_{2}$ demand lowers emissions in the $H_{2}$ sector and helps the electricity sector reduce its emissions. This is because flexible electrolysis helps to accommodate non-pollutive production into the electric system more efficiently.

\FloatBarrier

\begin{figure}[h]
    \centering
    \includegraphics[width=\FigureScale\textwidth]{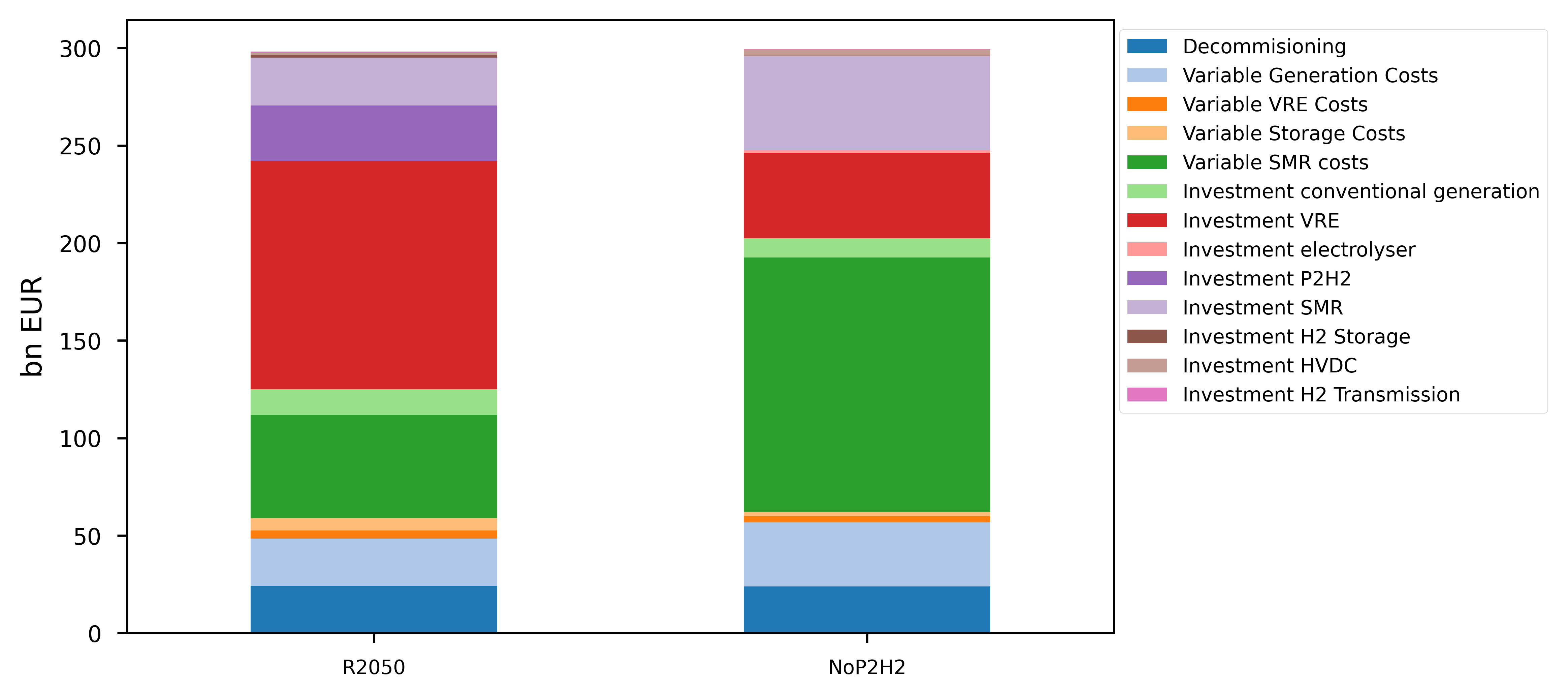}
    \caption{EU total system costs in R2050 and NoP2H2}
    \label{fig:EU_tot_syst_cost}
\end{figure}

Figure~\ref{fig:EU_tot_syst_cost} shows the total system cost for R2050 and NoP2H2. Interestingly, although the total costs between the two scenarios are similar (R2050 being 0.4\% lower), there is a significant redistribution of costs and R2050 yields 35\% lower emissions; see Figure \ref{fig:EU_H2_trans_invest}. As expected, R2050 shows significantly higher investments in P2H2 and VRE, while incurring significantly lower SMR (investment and variable) costs and variable generation costs.

\FloatBarrier
\subsection{Scenario variants}
\label{subsec:H2electr_ScnVar}

Table~\ref{tab:ref_sc_and_variants} presents different variants of the reference scenario R2050. By comparing these variants, we can separate the effect of other aspects on the system: variant NoP2H2 does not allow P2H2, thus showing the impact of electrifying $H_{2}$ demand, as widely discussed in the previous section, and some results are shown again here for the sake of completeness. The variant NoH2Storage variant allows us to observe the effect of $H_{2}$ flexibility through storage (i.e., time shifting), by not allowing investment in $H_{2}$ storage. Forcing the countries to only import/export energy via electricity in variant NoH2Transmission highlights the impact of $H_{2}$ transmission. The last scenario variant, NoETransmission, does not allow new expansions in the electricity network, thus only using the forecasted electricity transfer capacities; this previous variant analyses how investing in a $H_2$ network compares to expanding the current electricity network.

\begin{table}[h]
    \begin{tabularx}{\textwidth}{>{\raggedright\arraybackslash}X*{4}{>{\centering\arraybackslash}X}}
    \hline
    Scenario variants & Power-2-$H_{2}$ (electrolysis) & $H_{2}$ storage & $H_{2}$ transmission retrofit and new pipelines & Electrical transmission\footnotemark{} \\
    \hline
    R2050 & \checkmark & \checkmark & \checkmark & \checkmark \\
    NoP2H2 & X & \checkmark & \checkmark & \checkmark \\
    NoH2Storage & \checkmark & X & \checkmark & \checkmark \\
    NoH2Transmission & \checkmark & \checkmark & X & \checkmark \\
    NoETransmission & \checkmark & \checkmark & \checkmark & X \\
    \hline
    \end{tabularx}
    \caption{Scenario variants}
    \label{tab:ref_sc_and_variants}    
\end{table}

\footnotetext{Existing and forecasted initial electrical infrastructure is utilised from \citet{entsoe2020a}.}
\FloatBarrier

\subsubsection{$H_{2}$ supply and storage}
\label{subsubsec:H2electr_ScnVar_H2SuppStor}
\begin{figure}[h]
    \centering
    \includegraphics[width=\FigureScale\textwidth]{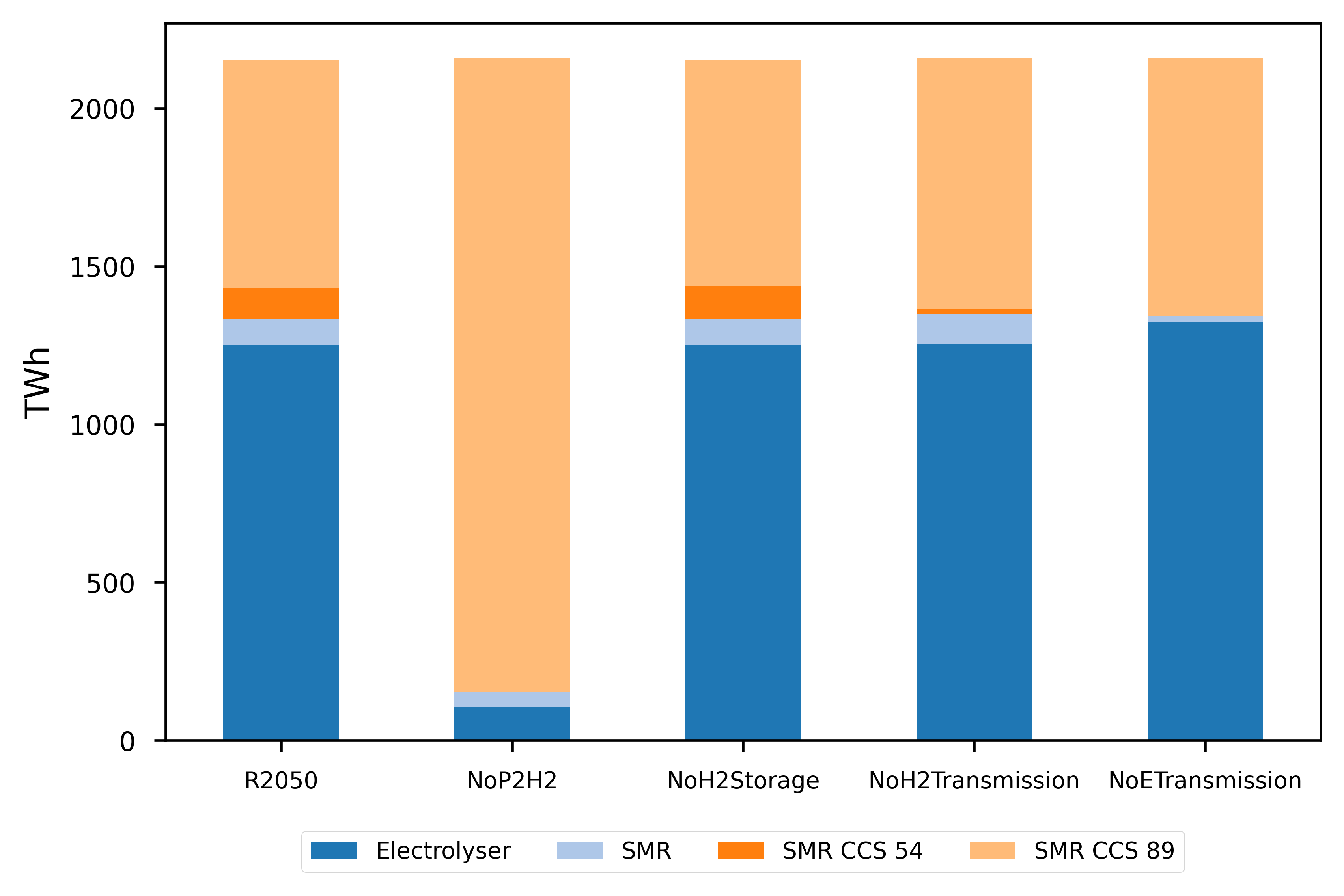}
    \caption{EU $H_2$ generation - comparison of R2050 and scenario variants}
    \label{fig:EU_H2_generation}
\end{figure}

Figure~\ref{fig:EU_H2_generation} presents the $H_2$ balances across various scenario variants at the EU level. The results indicate that in the R2050, NoH2Storage, NoH2Transmission and NoETransmission scenarios, the demand for electrified $H_2$ ranges from 58\% to 61\%. Notably, it is still optimal to satisfy approximately 40\% of the total $H_2$ demand through SMR technology with 89\% of carbon capture, even at a high $CO_{2}$ price of 250~\euro{}/ton. Furthermore, the policy limit of up to 50\% SMR generation is not reached at the EU level.

In the NoH2Storage variant, SMR with 54\% $CO_{2}$ capture appears viable in some countries. However, this option only meets 5\% of the total $H_2$ demand and the SMR technology with 89\% $CO_{2}$ capture is the preferred choice across all variants.

The NoH2Transmission variant exhibits the highest production of SMR without CCS, producing 97 TWh of $H_2$. Moreover, 80\% of this $H_2$ is produced by Germany, which is not allowed to invest in CCS technologies and cannot import $H_2$ from other countries.

\begin{figure}[h]
    \centering
    \includegraphics[width=\FigureScale\textwidth]{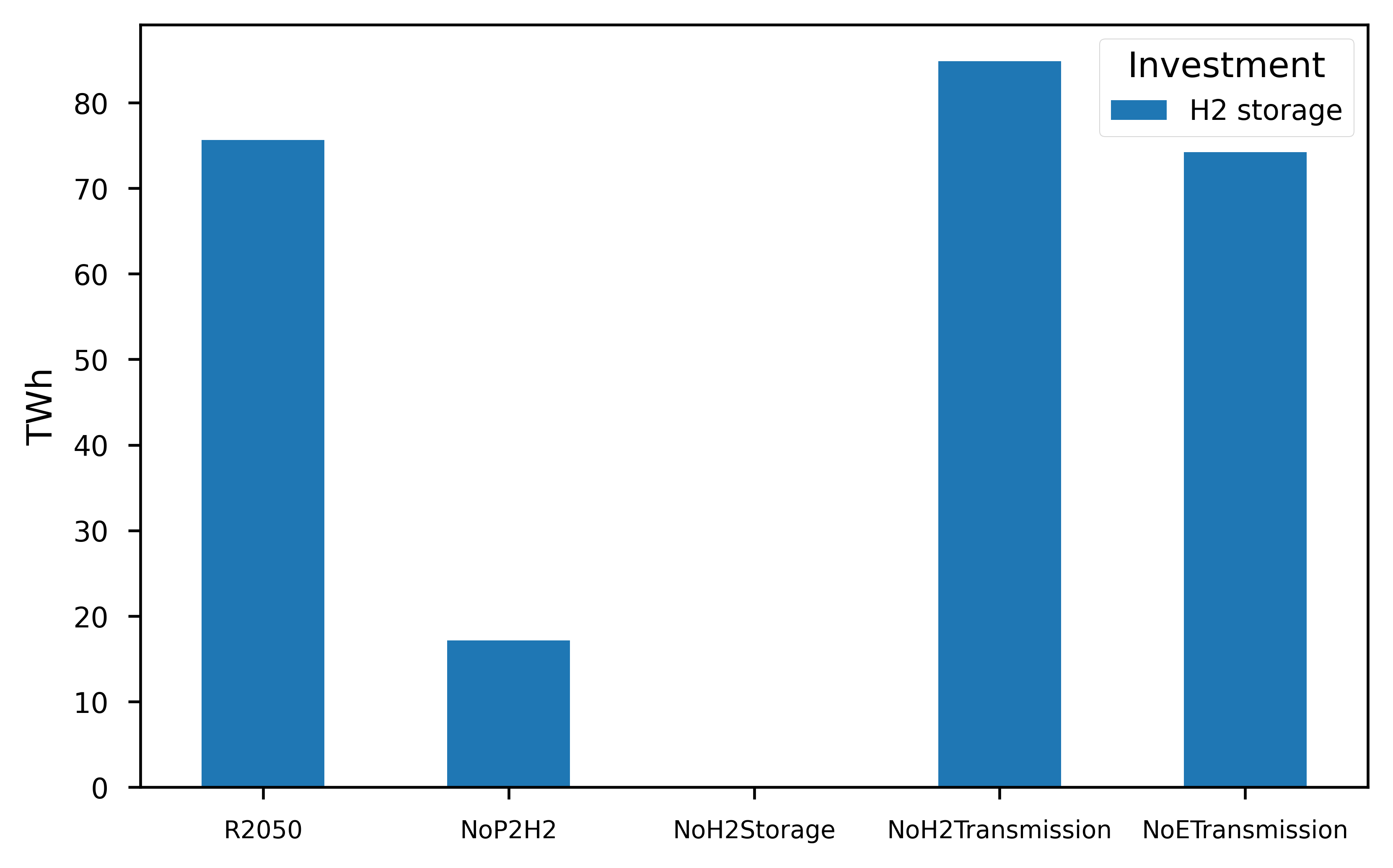}
    \caption{EU $H_2$ storage investments - comparison of R2050 and scenario variants}
    \label{fig:EU_H2_stor_invest}
\end{figure}

Figure~\ref{fig:EU_H2_stor_invest} displays the underground $H_2$ storage investments for the scenario variants relative to the reference scenario R2050. The following key observations can be made from the results. First, the underground $H_2$ storage requirement reduces significantly when there is no electrification of the $H_2$ demand. Specifically, in the NoP2H2 scenario, 58.4 TWh less storage is required compared to the R2050 scenario. The constant $H_2$ supply from SMR eliminates the need for time-shifting VRE-based $H_2$. Second, in the NoH2Transmission scenario, where $H_2$ transmission is not allowed, there is a 12\% increase (9.2 TWh) in $H_2$ storage requirements compared to the R2050 scenario. This increase is due to the higher need for time-shifting of $H_2$ to compensate for the flexibility lost from geographical shifting.

\FloatBarrier
\subsubsection{Electricity supply and demand}
\label{subsubsec:H2electr_ScnVar_PowerDemSup}
\begin{figure}[h]
    \centering
    \includegraphics[width=\FigureScale\textwidth]{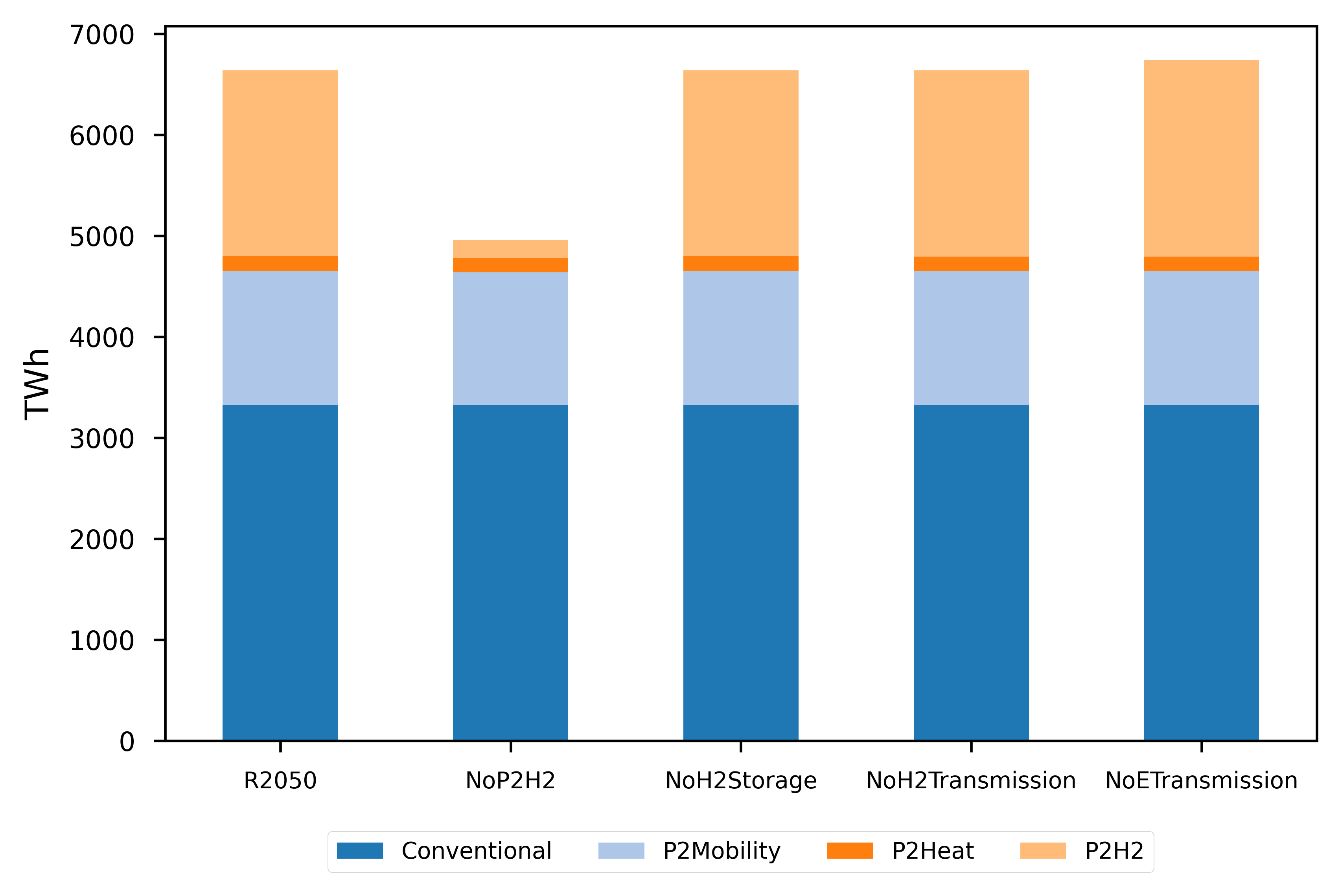}
    \caption{EU electricity} demand - comparison of R2050 and scenario variants
    \label{fig:EU_power_demand_sc_var}
\end{figure}

\begin{figure}[h]
    \centering
    \includegraphics[width=\FigureScale\textwidth]{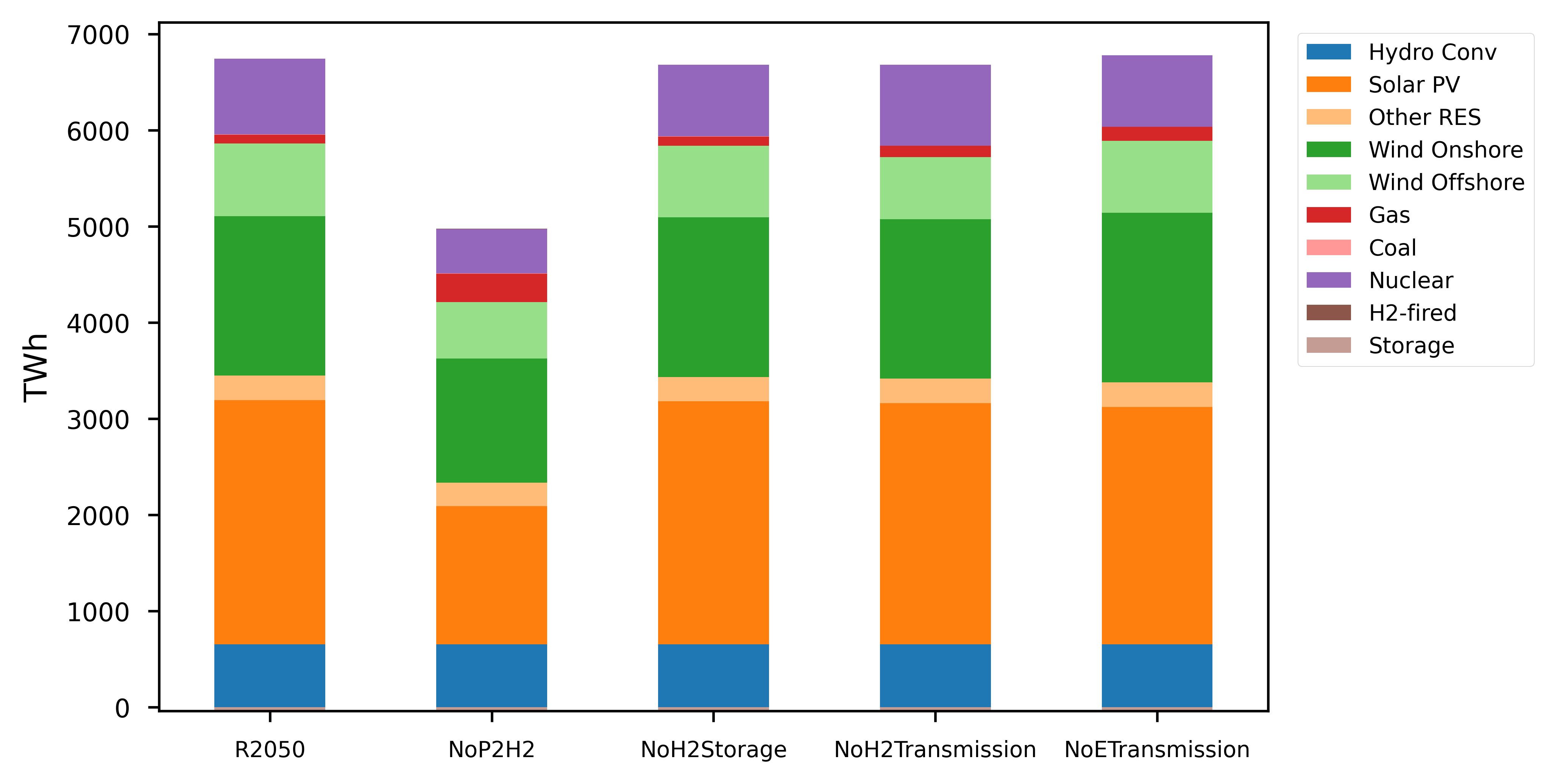}
    \caption{EU electricity} generation - comparison of R2050 and scenario variants
    \label{fig:EU_power_generation_sc_var}
\end{figure}

Figure~\ref{fig:EU_power_demand_sc_var} and Figure~\ref{fig:EU_power_generation_sc_var} compare the EU electricity demand and supply of the scenario variants to the reference scenario R2050. The results show that the $H_{2}$ electrical demand is similar for the scenarios R2050, NoH2Storage, and NoH2Transmission, with only a 3 TWh difference. This marginal demand change can be attributed to the system's ability to use spatial or temporal flexibility to achieve similar levels of $H_{2}$ electrification. In addition, the $H_{2}$ electrification levels and energy mixes are similar for the R2050 and NoETransmission scenarios, suggesting that a system with the expected transmission expansion by 2050 is already near the optimal solution. Finally, in comparison to R2050, the $H_{2}$ electrification levels are similar in the NoH2Storage and NoH2Transmission scenarios (all around 6638 TWh); however, VRE production drops 0.5\% (21 TWh) and 3\% (140 TWh), respectively. This highlights how the ability of $H_{2}$ to follow VRE production, either in time or space, can help to increase VRE production by raising flexible $H_{2}$ electrification levels. In contrast, nuclear production increases by 5\% (40 TWh) in NoH2Storage, indicating that nuclear energy may be a better alternative than VRE for electrifying more inflexible $H_{2}$. The differences in VRE and nuclear production affect the total system costs. For instance, Section \ref{subsubsec:H2electr_ScnVar_TotSysCost} shows that the NoH2Storage and NoH2Transmission variants are around 10 bn\euro{} and 4 bn\euro{} more expensive than the reference scenario.

\FloatBarrier
\subsubsection{Energy transmission}
\label{subsubsec:H2electr_ScnVar_EnergyTransm}
\begin{figure}[h]
    \centering
    \includegraphics[width=\FigureScale\textwidth]{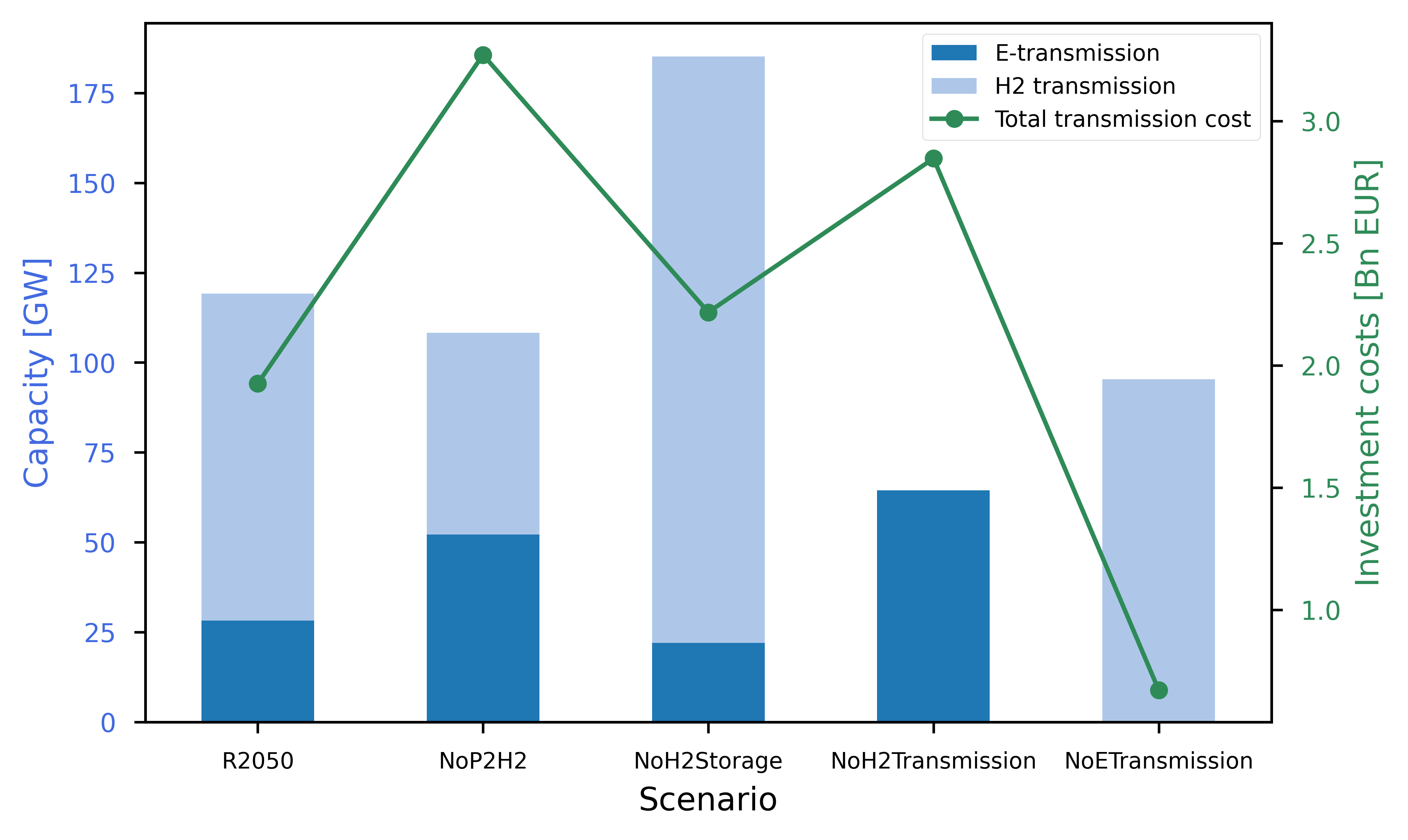}
    \caption{EU transmission investments and costs comparison of R2050 and scenario variants}
    \label{fig:EU_transm_invest_sc_var}
\end{figure}

Figure~\ref{fig:EU_transm_invest_sc_var} presents the required electrical and $H_2$ transmission results for the reference scenario (R2050) and its variants, along with their corresponding investment costs. Section \ref{sec:opt_model} explains that COMPETES calculates the optimal transmission infrastructure needed to couple the demand and supply of electricity and $H_2$ among different countries.

In the R2050 case, the required $H_2$ transmission decreases significantly from 90 GW in the NoP2H2 scenario to 56 GW, which is expected since the model assumes unlimited SMR potential within each country. However, the total costs in this variant are the highest due to the required expansion of the electricity network. This decrease in $H_2$ transmission requirements is also driven by the lack of green $H_2$ production from countries with VRE resources, such as Spain and France.

In the NoH2storage scenario, the required expansion on the electricity network decreases by 6 GW compared to the R2050 case, while $H_2$ transmission increases by 72 GW. This reduces total costs compared to the R2050 case, as an extra investment in the electricity network is avoided, which is more costly (per GW) than expanding the $H_2$ network. Not allowing $H_{2}$ storage also results in increasing $H_{2}$ transmission capacity by 80\%, which uses pipelines to store $H_{2}$ to be used in high-demand moments.

In the NoETransmission scenario, which only uses the future forecast of transfer capacities, the results show that only 5GW of extra investment in $H_2$ transmission is required. The associated energy transmission investment costs are reduced by 180\% compared to the R2050 case.

In the NoH2Transmission scenario, where there is no possibility to invest in $H_2$ transmission and the expansion of the electricity network is driven naturally, resulting in an increase in electrical network investments by 130\% compared to the R2050 case. As shown in Figure 17, this expansion in the electricity network is complemented by higher investment in $H_2$ storage.

\FloatBarrier
\subsubsection{$CO_{2}$ emissions}
\label{subsubsec:H2electr_ScnVar_CO2Emiss}
\begin{figure}[h]
    \centering
    \includegraphics[width=\FigureScale\textwidth]{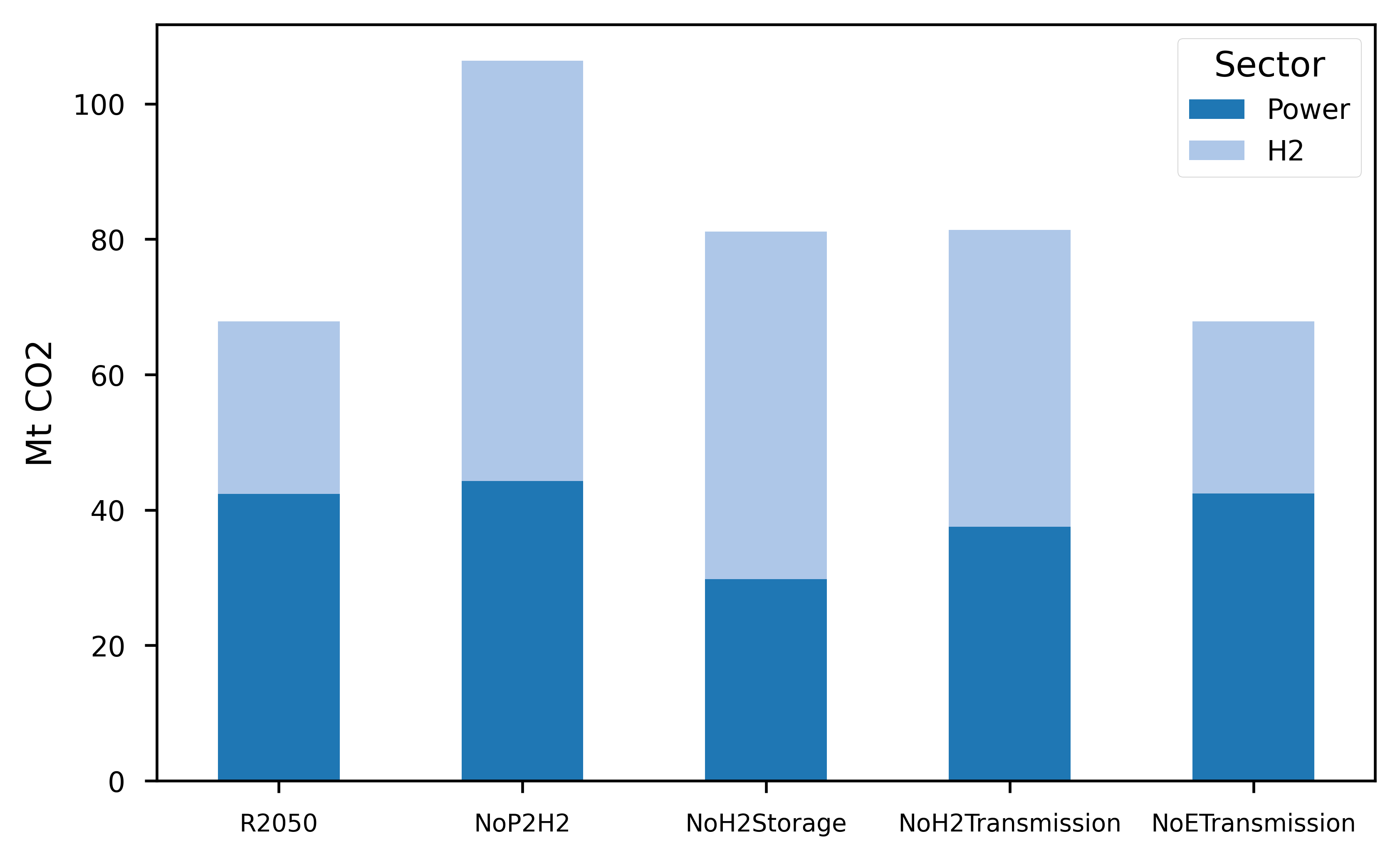}
    \caption{EU $CO_{2}$  emissions - comparison of R2050 and scenario variants}
    \label{fig:EU_CO2_emiss_sc_var}
\end{figure}

The results of the scenario variants regarding the $CO_{2}$ emissions in the EU electricity and $H_2$ sectors are presented in Figure~\ref{fig:EU_CO2_emiss_sc_var}. The total $CO_{2}$ emissions increased by 19.5\% and 20\% in the NoH2Storage and NoH2Transmission cases, respectively, compared to the R2050 case. The increase in $CO_{2}$ emissions in the $H_{2}$ sector can be attributed to the rise in the use of SMR technologies to supply the $H_{2}$ demand in both cases.

Interestingly, $CO_{2}$ emissions are similar in the R2050 and the NoETransmission case. This suggests that further expanding the electricity network does not necessarily result in $CO_{2}$ emissions reductions after achieving a specific transfer capacity between countries.

\FloatBarrier
\subsubsection{Total system costs}
\label{subsubsec:H2electr_ScnVar_TotSysCost}
\begin{figure}[h]
    \centering
    \includegraphics[width=\FigureScale\textwidth]{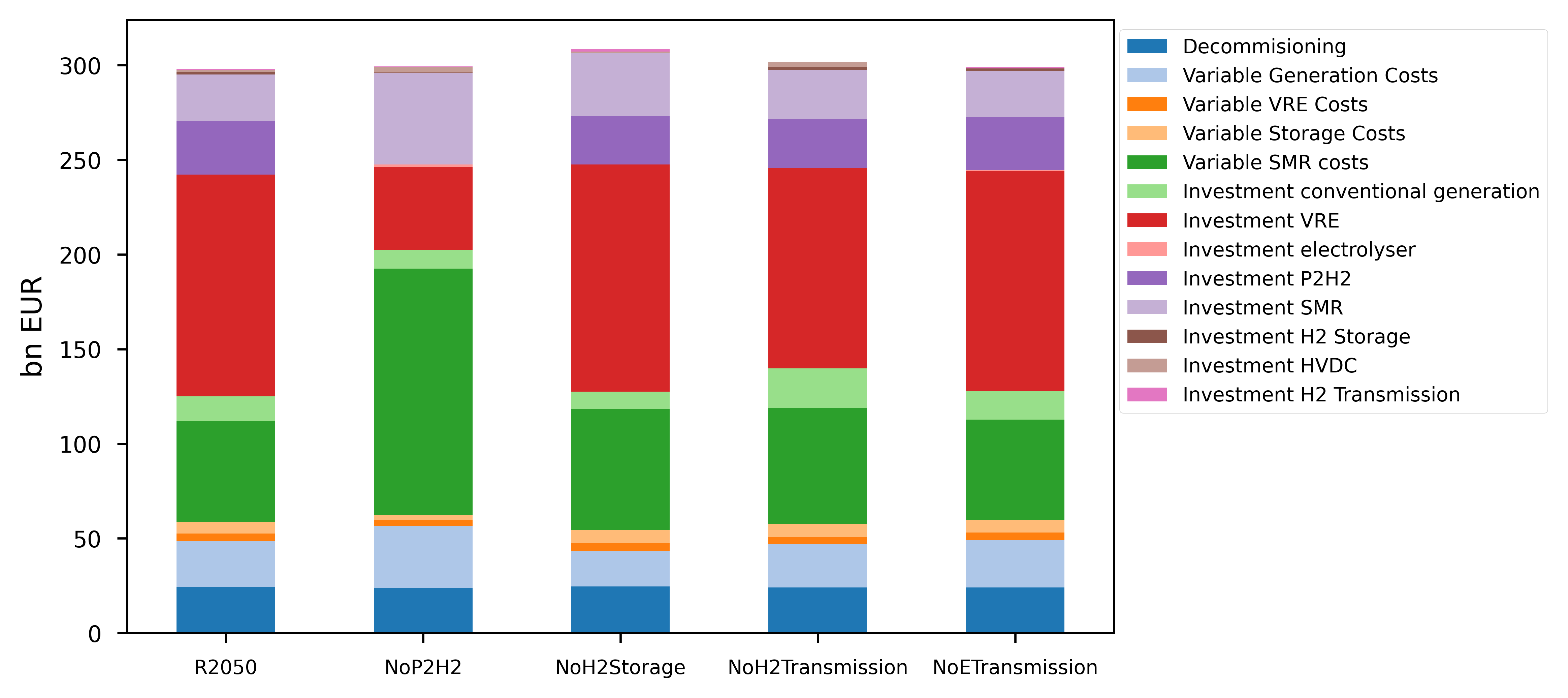}
    \caption{EU total system costs - comparison of R2050 and scenario variants}
    \label{fig:EU_tot_sys_cost_sc_var}
\end{figure}

Figure~\ref{fig:EU_tot_sys_cost_sc_var} shows the system cost distributions for the R2050 scenario and variants. The R2050 scenario has the lowest total system costs, as it can access all investment options to achieve an optimal solution. It can minimise system costs by using the optimal combination of electrification, storage, and transmission. In contrast, the NoH2Storage case has the highest total system costs, about 3.4\% (10.3 b\euro{}) higher than in the R2050 scenario. This is because the NewVRE costs, i.e., investment in wind and solar, are lower, while the variable $H_{2}$ costs increase due to the gas costs of $H_2$ production via SMR.

\begin{figure}[h]
    \centering
    \includegraphics[width=\FigureScale\textwidth]{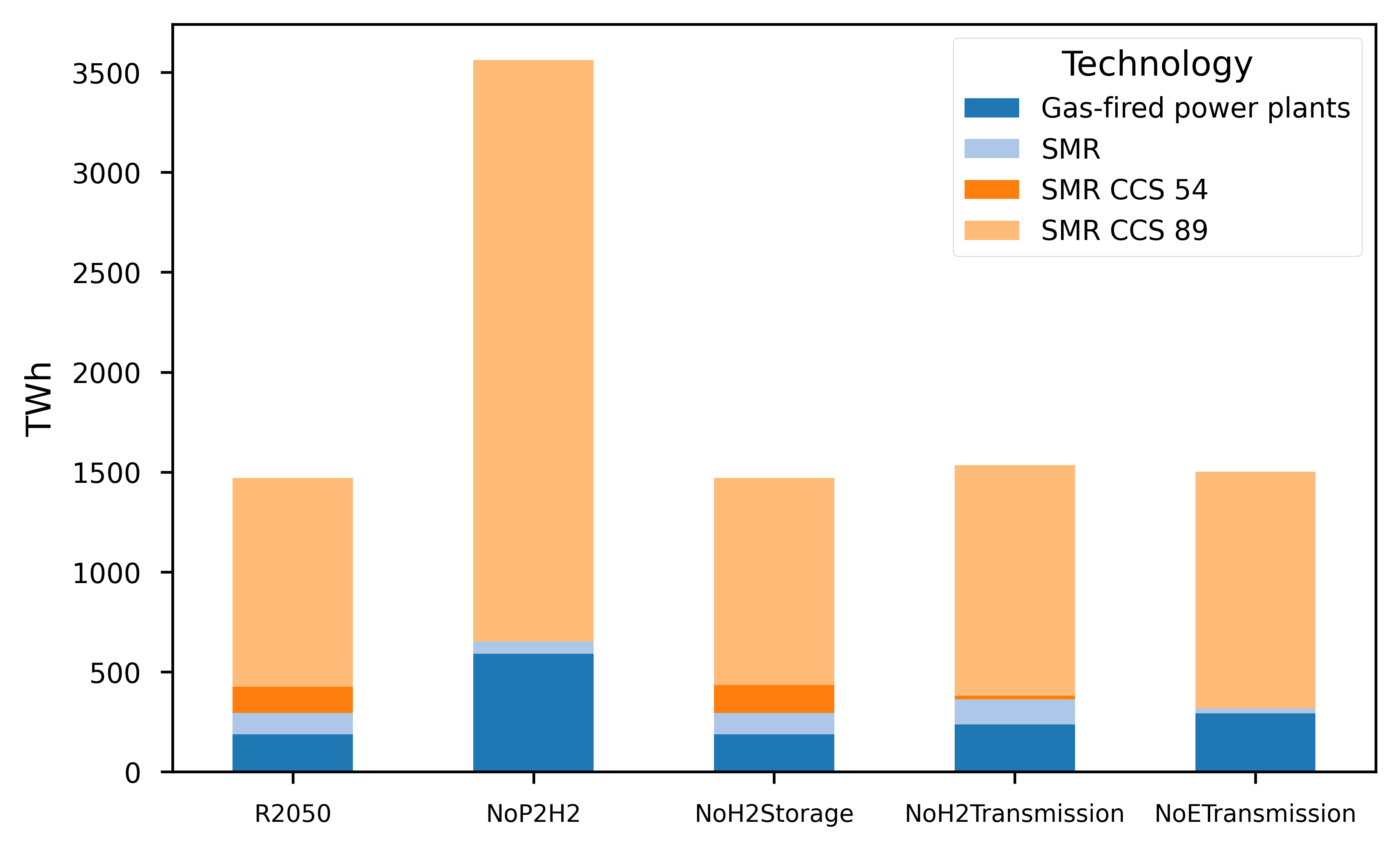}
    \caption{EU total gas consumption - comparison of R2050 and scenario variants}
    \label{fig:EU_tot_gas_consumption_sc_var}
\end{figure}

In the NoP2H2 scenario, costs are significantly shifted from investment in VRE to variable SMR costs. This is due to the lack of $H_{2}$ electrification, thus requiring less VRE capacity and higher gas consumption to produce $H_{2}$. Figure \ref{fig:EU_tot_gas_consumption_sc_var} shows that the gas consumption in the NoP2H2 is around 2.4 times higher than the other scenarios. There is also an increase in gas-fired power plant output, resulting in a rise of 8.5 b\euro{} in variable generation costs (fuel costs) shown in Figure~\ref{fig:EU_tot_sys_cost_sc_var}. Moreover, the NoETransmission scenario shows only a 0.3\% total system cost increase compared to the R2050 scenario. In this case, there is a decrease in the variable VRE costs, but an increase in the conventional generation costs. Since the electricity network cannot expand further, it is not possible to integrate more VRE, so traditional generation is required.

\FloatBarrier
\subsubsection{Electricity vs $H_{2}$ Transport}
\label{subsubsec:H2electr_ScnVar_ElecVsH2Transp}
\begin{figure}[h]
    \centering
    \includegraphics[width=\FigureScale\textwidth]{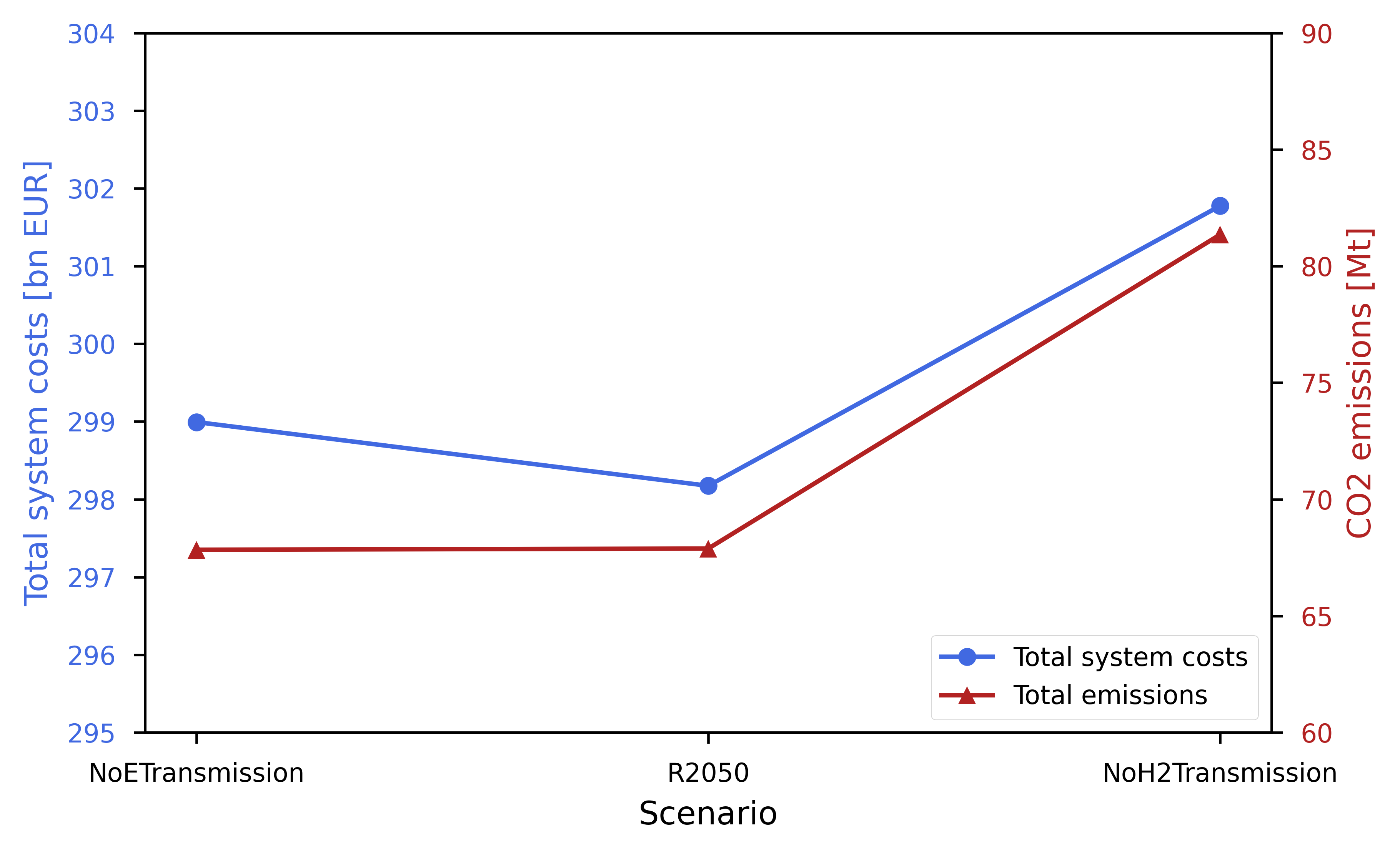}
    \caption{$H_{2}$ vs. Electricity transmission}
    \label{fig:H2_vs_elec_transm_sc_var}
\end{figure}

In this section, we compare $H_2$ and electricity transport in terms of total system costs and $CO_{2}$ emissions; see Figure~\ref{fig:H2_vs_elec_transm_sc_var}. To this end, we analyse three transmission scenario variants: the R2050 reference scenario, the NoETransmission, and the NoH2Transmission extreme cases. The R2050 scenario represents the optimal tradeoff between electricity and $H_{2}$ transmission expansion; the NoETransmission scenario assumes that no further development of the electricity network is allowed; and the NoH2Transmission variant considers the absence of $H_2$ transport via pipelines.

Our findings reveal that relying solely on electricity transmission, as in the NoH2Transmission scenario, leads to the worst outcome, resulting in higher costs and $CO_{2}$ emissions. Specifically, NoH2Transmission increases costs by 1.2\% and emits around 20\% more $CO_{2}$ compared to R2050. Conversely, the NoETransmission scenario has comparable $CO_{2}$ emissions to the optimal mix (R2050), but system costs are slightly higher (0.3\%, 0.8 b\euro{}). Notably, investing solely in the $H_{2}$ network and storage eliminates the need for an additional 28 GW of electrical network expansion, as shown in Figure~\ref{fig:EU_transm_invest_sc_var}. Our results indicate that a system with the expected transmission expansion by 2050 is already very close to the optimal solution. Therefore, to achieve more significant cost savings and $CO_{2}$ emissions reductions, we suggest policymakers should focus on facilitating $H_{2}$ transport through pipelines rather than increasing electricity transport.

\FloatBarrier
\subsection{Revisiting the research questions}
This section has provided solutions to the research questions in Section \ref{subsec:contribution}. We summarise our findings as follows:

\begin{itemize}
    \item \textit{Is it possible to effectively include the retrofit of natural gas networks as an investment option in energy planning optimisation models? Yes--} The option of retrofitting natural gas networks as an investment option has been considered in COMPETES, according to the proposed formulation in Section \ref{sec:retrofit_formulation}. The formulation's effectiveness lies in its LP mathematical nature, allowing for solving large-scale optimisation models, like the European reference scenario and its variations, without a significant increase in the computational time burden.
    \item \textit{What are the impacts on the total costs and $CO_2$ emissions of retrofitting the existing gas infrastructure for $H_2$ transport in the EU by 2050?} We can measure the impact by comparing the reference scenario (R2050) with its variants, particularly the scenario where no $H_{2}$ investment is made in retrofitting or new pipelines (NoH2Transmission). This results in a 1.2\% increase in costs compared to R2050 and around 20\% more $CO_2$ emissions.
    \item \textit{How do investments in new electrolysers, hydrogen transmission, and storage infrastructure impact total $CO_2$ emissions compared to scenarios where these investments are not made?} We have created specific scenario variants to answer particular questions. For example, the NoP2H2 variant shows us the impact of not investing in new electrolysers. The NoH2Transmission variant shows the effect of not retrofitting or building new pipelines for $H_2$ transport, while the NoH2Storage variant focuses on $H_2$ storage. Lastly, the NoETransmission variant measures the impact of not having additional electricity transmission. To quantify the impacts, we compared the results of each variant to the reference scenario. For instance,  the total $CO_2$ emissions increased by around 20\% the NoH2Storage and NoH2Transmission cases compared to the R2050 case, highlighting the importance of the sector coupling between the power and hydrogen sectors to lower the total $CO_2$ emissions.
\end{itemize}

\FloatBarrier
%---------------------------------------------------------------------------
%           Conclusions
%---------------------------------------------------------------------------
\section{Conclusions}
\label{sec:conclusions} 
This paper has examined the potential of hydrogen ($H_{2}$) electrification to transform the power systems of the European Union in the year 2050. Various scenarios were used to measure the impact of essential investment decisions, such as electrolysers,  storage, retrofitting, and new transmission systems. One key finding in our research is that a strategic balance of $H_{2}$ electrification and Steam Methane Reforming (SMR), alongside efficient transmission of both electricity and $H_{2}$, helps to reduce $CO_{2}$ emissions and enables a sustainable and cost-effective power system in the EU for 2050. Moreover, electrolysers might become essential tools in the power sector, providing flexibility to replace peak units like gas with non-polluting technologies. We, therefore, conclude that electrifying the $H_{2}$ production is essential to meet the European Union's goals of reducing long-term emissions.

The proposed formulation to consider retrofitting natural gas networks as an investment opportunity in COMPETES has offered valuable information. Without the option to transport $H_{2}$ through either retrofitting or new pipelines, there is a significant increase in costs and $CO_{2}$ emissions. Finding a balance between transporting $H_{2}$ and electricity is crucial to create an optimal system. The expansion of electricity transmission by 2050 is already near the optimal solution, so the focus should be on facilitating $H_2$ transport, including retrofitting.

In conclusion, the results in this paper provide a measurable assessment that can guide future policy decisions on this issue, and it is a valuable resource for policymakers to make informed decisions.

\section{Acknowledgements}
This research has received funding from the Fuel Cells and Hydrogen 2 Joint Undertaking under grant agreement No. 735503 (H2Future project) and from the European Climate, Infrastructure and Environment Executive Agency (CINEA) under the European Union's HORIZON Research and Innovation Actions under grant agreement No. 101095998.

The work of Mr. Morales-España was partially funded by the European Climate, Infrastructure and Environment Executive Agency under the European Union’s HORIZON Research and Innovation Actions under grant agreement N°101095998.

\textbf{Disclaimer}: Views and opinions expressed are, however, those of the author(s) only and those of the European Union or CINEA. Neither the European Union nor the granting authority can be held responsible for them.

Finally, the authors express their sincere gratitude to Lauren Clisby for her review and editing of this manuscript.

\section{Author Contributions}
\textbf{Germán Morales-España}: Methodology, Software. \textbf{Ricardo Hernández-Serna}: Data Curation, Writing, Investigation. \textbf{Diego A. Tejada-Arango}: Visualization, Writing, Reviewing, and Editing. \textbf{Marcel Weeda}: Supervision, Project administration.

%% The Appendices part is started with the command \appendix;
%% appendix sections are then done as normal sections
\appendix

%% If you have a bibdatabase file and want BibTeX to generate the
%% bibitems, please use
%%
\bibliographystyle{elsarticle-num-names} 
\bibliography{references}

%% else use the following coding to input the bibitems directly in the
%% TeX file.

% \begin{thebibliography}{00}

% %% \bibitem[Author(year)]{label}
% %% Text of bibliographic item

% \bibitem[ ()]{}

% \end{thebibliography}
\end{document}